\documentclass[11pt]{article}

\usepackage{amssymb}
\usepackage{epsfig}

\newtheorem{Lemma}{Lemma}
\newtheorem{Proposition}[Lemma]{Proposition}

\newtheorem{Theorem}[Lemma]{Theorem}

\newtheorem{Corollary}[Lemma]{Corollary}

\newcommand{\Aopt}{A^{\mbox{{\tiny opt}}}}
\newcommand{\Bopt}{B^{\mbox{{\tiny opt}}}}
\newcommand{\Anopt}{A_n^{\mbox{{\tiny opt}}}}
\newcommand{\Bnopt}{B_n^{\mbox{{\tiny opt}}}}
\newcommand{\Atopt}{A_{t-1}^{\mbox{{\tiny opt}}}}
\newcommand{\Btopt}{B_{t-1}^{\mbox{{\tiny opt}}}}
\newcommand{\ATopt}{A_{T-1}^{\mbox{{\tiny opt}}}}
\newcommand{\BTopt}{B_{T-1}^{\mbox{{\tiny opt}}}}
\newcommand{\med}{\mbox{median}}

\newcommand{\ed}{\ \stackrel{d}{=} \ }
\newcommand{\cd}{\ \stackrel{d}{\rightarrow} \ }

\renewcommand{\SS}{\mbox{${\cal S}$}}

\newcommand{\eps}{\varepsilon}

\newcommand{\bx}{{\mathbf x}}
\newcommand{\by}{\mathbf y}
\newcommand{\bb}{\mathbf b}
\newcommand{\Z}{{\mathbb Z}}

\newcommand{\Reals}{{\mathbb R}}
\newcommand{\R}{{\mathbb R}}

\newcommand{\Bcal}{\mathcal B}
\newcommand{\Vcal}{{\mathcal V}}
\newcommand{\Wcal}{{\mathcal W}}

\newcommand{\bV}{\widetilde{V}}
\newcommand{\bW}{\widetilde{W}}
\newcommand{\bM}{\widetilde{M}}

\newcommand{\sfrac}[2]{{\textstyle\frac{#1}{#2}}}
\newcommand{\qed}{\ \ \rule{1ex}{1ex}}

\def\ux{\underline{x}}

\newcommand{\Pb}{{\mathrm{P}}}
\newcommand{\Eb}{{\mathrm{E}}}
\newcommand{\Zbold}{{\mathbb{Z}}}

\def\ind{{\rm 1\hspace{-0.90ex}1}}
\def\oG{\overline{G}}
\def\oF{\overline{F}}

\newcommand{\Ical}{\mathcal{I}}

\author{David J. Aldous\thanks{Research supported by N.S.F. Grant
DMS0704159}
\\
\\
       University of California\\
       Department of Statistics\\
        367 Evans Hall \# 3860\\
       Berkeley CA 94720-3860\\
       aldous@stat.berkeley.edu
\and Charles Bordenave\thanks{Research supported by N.S.F. Grant CCF-050023.}
 \\
\\
       University of California\\
       Departments of EECS and Statistics\\
        257 Cory Hall  \# 1770\\
       Berkeley CA 94720-1770\\
       charles.bordenave@berkeley.edu
\and Marc Lelarge\\
\\
       INRIA-ENS\\
       45 rue d'Ulm\\
       75230 Paris Cedex 5\\
       marc.lelarge@ens.fr}

\title{Dynamic Programming Optimization over Random Data: 
the Scaling Exponent for Near-optimal Solutions}

\begin{document}

\maketitle

\begin{abstract}
A very simple example of an algorithmic problem solvable by dynamic programming is to maximize,
over $A \subseteq \{1,2,\ldots,n\}$, the objective function 
$|A| - \sum_i \xi_i \ind(i \in A,i+1 \in A)$ 
for given $\xi_i > 0$.  
This problem, with random $(\xi_i)$, provides a test example for studying the relationship between optimal 
and near-optimal solutions of combinatorial optimization problems.  
We show that, amongst solutions differing from the optimal solution in a small proportion $\delta$ 
of places, we can find near-optimal solutions whose objective function value differs from the optimum 
by a factor of order $\delta^2$ but not smaller order. 
We conjecture this relationship holds widely in the context of dynamic programming over random data,
and Monte Carlo simulations for the Kauffman-Levin NK model  
are consistent with the conjecture.
This work is a technical contribution to a broad program initiated in Aldous-Percus (2003) 
of relating such scaling exponents to the algorithmic difficulty of optimization problems.

\end{abstract}

{\bf Key words.}  
Dynamic programming, 
local weak convergence,
Markov chain, 
near-optimal solutions,
optimization,
probabilistic analysis of algorithms,
scaling exponent.

{\bf AMS subject classifications.}  68Q25,  90C39, 60J05.

\newpage
\section{Introduction and Motivation}

\subsection{Near-optimal solutions in combinatorial optimization}

Consider a 
combinatorial optimization problem 
which is ``size $n$" 
in the sense that a feasible solution 
$\bx = (x_i, 1 \leq i \leq n)$
consists of $n$ elements 
(e.g. edges of a graph; binary digits)
subject to some constraints, 
and the objective function $f(\bx)$ 
is akin to a sum over $i$ of costs or rewards associated with each $x_i$.
In such a setting one can define the relative distance between the structure of a feasible solution 
$\bx$ and the optimal solution $\bx^*$ by
\[
\delta_n({\bf x}) = % \mbox{\{number of edges in {\bf x}
n^{-1} |\{i: x_i \neq x^*_i\}|
\]
and the relative difference
in objective function is
$n^{-1}|f(\bx) - f(\bx^*)|$.
So the quantity
\begin{equation}
\eps_n(\delta) := 
\min 
\{n^{-1}|f(\bx) - f(\bx^*)|
\ : \delta_n(\bx) \geq \delta \} 
\label{fmin}
\end{equation}
measures how close we can get to the optimal value using feasible solutions which have
non-negligibly different structure from the optimal solution.
A program initiated in
\cite{me103}
is to study this quantity for 
combinatorial optimization problems 
over {\em random} data.  
In this setting $\eps_n(\delta)$ becomes a random variable, 
but in many cases one expects
that as $n \to \infty$ there is a 
{\em deterministic} limit function $\eps(\delta)$.
Motivation for this program is a conjecture that 
(within some suitable class of problems)
\[ \eps(\delta) \asymp \delta^\alpha \mbox{ as } \delta \to 0 \]
for some {\em scaling exponent} $\alpha $,
whose value is 
robust under model details, and that 
for ``algorithmically easy" problems 
we have $\alpha = 2$
(which of course mimics the behavior we expect by calculus for smooth functions 
$f: \Reals^d \to \Reals$)
whereas
for ``algorithmically hard" problems 
we have $\alpha > 2$.
Here is the previous evidence in support of this conjecture.

(i) 
{\em Traveling salesman problem}
and
{\em minimum matching problem} 
\cite{me103}.
In the random link (mean-field) model, a cavity method analysis 
(non-rigorous but generally regarded as accurate) enables one to compute 
$\eps(\delta)$ numerically and to observe scaling exponent 
$\alpha = 3$.
In the random Euclidean model, Monte Carlo simulations suggest the same 
$\alpha = 3$.

(ii) 
{\em Minimum spanning tree}.
Here we expect $\alpha = 2$.
This is proved in \cite{me115} 
for the $d \geq 2$ dimensional
random Euclidean model and also for a ``disordered lattice" model.

The purpose of this paper is to consider some problems which are 
algorithmically easy to solve via dynamic programming, 
and where we therefore expect $\alpha = 2$.
We first give a trivial but instructive case (section \ref{sec-2})
and then describe a prototypical ``interesting" case, the 
Kauffman-Levin $NK$ model (section \ref{sec-NK}).
Here both a heuristic argument and simulations suggest $\alpha = 2$, 
but we do not have a proof.
Our main focus is on giving a complete analysis of a simple non-trivial model 
(section \ref{sec-4})
where we are required to pick a subset 
$A \subseteq 
[n] :=
\{1,2,\ldots,n\}
$ of items with a reward of $1$ per item picked and i.i.d. costs 
$\xi_i$ incurred if both items $i$ and $i+1$ are picked. 
Theorem \ref{th2} establishes $\alpha = 2$ for this specific model. 
In these dynamic programming examples and the minimum spanning tree example, 
the key structural property is that the near-optimal solutions attaining the minimum in (\ref{fmin}) 
differ from the optimal solution via only ``local changes", each local change affecting only a number 
of items which remains $O(1)$ as $\delta \to 0$.  
It is natural to speculate that this structural property corresponds quite generally to the $\alpha = 2$ case.

\paragraph{Related work} 

We do not know any other lines of research in theoretical computer science which are close to 
the topic of this paper.  
A recent survey of average-case complexity of NP problems is given in \cite{0606037}.  
Interest in the average-case  gap between optimal and second-optimal solutions arises in several contexts, e.g.  \cite{1007409}.  
Closer in spirit is the statistical physics of disordered systems, where for low temperatures 
the Gibbs distribution on configurations concentrates on near-minimal-cost configurations.  
%In this context it is natural to consider the {\em overlap} between two random picks from the Gibbs distribution. 
In the context of random energy models (the precise analog of optimization over random data), two random picks from the Gibbs distribution over the same random choice of energy are called {\em replicas}, and study of such replicas and their overlaps is a central theme of the {\em replica method} \cite{MPV87,MR1993891}.  
So that topic studies the structural difference between two {\em typical} near-optimal configurations, whereas we study the {\em maximal} (over all near-optimal configurations) structural difference from the optimal configuration.  
Our mathematical arguments are much less sophisticated than those in statistical physics, but there are some intriguing parallels, described briefly in section \ref{sec-cavity}.

\subsection{A trivial example}
\label{sec-2}
Let $(X_i, i \geq 1)$ be i.i.d. real-valued random variables 
with continuous density $h(x)$ and $\Eb X < \infty$.
For each $n$ consider the problem of finding
\[ M_n = \max_{A \subseteq [n]}
\sum_{i \in A} (X_i - 1) . \] 
The maximum is 
obviously obtained by choosing
$A = \{i: X_i > 1\}$ 
and then as $n \to \infty$
\[ n^{-1} M_n \to \Eb (X_1 - 1)^+ 
\mbox{ a.s. } \]
Fix $0 < \delta < 1$.
It is also obvious that the subset $A^\prime$ 
that minimizes
\[ M^\prime_n = \max_{A^\prime \subseteq [n]}
\sum_{i \in A^\prime} (X_i - 1)  \]
\[ \mbox{subject to $|A^\prime \bigtriangleup A| \geq \delta n$}\]
is the subset $A^\prime = A \bigtriangleup D$
where $D$ is the set of indices of the $\lceil \delta n \rceil$
smallest values of $|X_i - 1|$.
So as $n \to \infty$
\[ n^{-1}(M_n - M^\prime_n) \to_{L_1}
\eps(\delta) := 
\int_{1-a(\delta)}^{1+a(\delta)}
|x-1| h(x) \ dx \]
where $a(\delta)$ is defined by 
\[ \delta = 
\int_{1-a(\delta)}^{1+a(\delta)}
h(x) \ dx . \]
So by continuity of $h(x)$, and assuming 
$0<h(1)<\infty$,
as $\delta \downarrow 0$ we have
\begin{equation}
a(\delta) \sim \sfrac{\delta}{2h(1)}; \quad 
\eps(\delta) \sim a^2(\delta) h(1) \sim \sfrac{\delta^2}{4h(1)}
\label{iid-case}
\end{equation}
which is the desired ``scaling exponent $=2$" result.

\paragraph{Discussion.}
(i) This example illustrates a feature that arises in other examples, that proving $\alpha = 2$ reduces to
showing that the density of a certain measure at a certain point is finite and non-zero.  
In nontrivial examples the measure in question arises in the {\em analysis} of the problem
rather than the statement of the problem: see
Lemma \ref{le:mu} below 
and Proposition 8 of \cite{me115}.

\noindent
(ii) In this example we could see the form of the best near-optimal solution by inspection, but 
a systematic method is to use Lagrange multipliers.  
In this example, introduce a parameter $\theta > 0$ and consider for each $n$
 \[ A_\theta := \arg \max_A 
 \left( \sum_{i \in A} (X_i - 1) + \theta |A \bigtriangleup A^*| \right) \]
 where $A^* = \{i: X_i > 1\}$ is the optimal solution. 
 By inspection the solution is
 \[ A_\theta = \{i: \ 1-\theta \leq X_i \leq 1 \mbox{ or } 1+\theta \leq X_i\}  . \]
 Although now $|A_\theta \bigtriangleup A^*|$ is random, we can use the law of large numbers
 to obtain existence of the limits
 \begin{eqnarray*}
 \delta(\theta) &:=& \lim_{n \to \infty}n^{-1} |A^* \bigtriangleup A_\theta| = \int_{1- \theta}^{1+\theta}
h(x) \ dx  \\
\eps(\theta) &:=& \lim_{n \to \infty} n^{-1} \left( \sum_{i \in A^* } (X_i - 1) - \sum_{i \in A_\theta } (X_i - 1)\right) 
= \int_{1- \theta}^{1+\theta}
|x-1| \ dx  .
\end{eqnarray*}
By the interpretation of Lagrange multipliers, this is an implicit function representation of $\eps$ as a function of $\delta$, and rederives 
 the limit (\ref{iid-case}) above.

\subsection{The NK model}
\label{sec-NK}

The Kauffman-Levin NK model of random fitness landscape has attracted extensive literature 
in statistical physics  \cite{kauffman93,weinberger91}.
For our version of the model 
we fix $K \geq 2$.
We seek to minimize, over binary sequences
$\bx = (x_1,\ldots,x_N)$,
the objective function
$H_N(\bx) = \sum_{i=1}^{N-K} W_i(x_{i},x_{i+1},\ldots,x_{i+K})$,
where the values
$(W_i(b_0,b_1,\ldots,b_K):i \geq 1, \bb \in \{0,1\}^{K+1})$
are independent exponential($1$) random variables.  
This is algorithmically easy via dynamic programming.
Write $\bx^N$ for the minimizing sequence.
By subadditivity
there is an a.s. limit
$N^{-1} H_N(\bx^N) \to c_K$.
For a general sequence $\by = \by^N$ write
\begin{eqnarray*}
\delta_N(\by) &=& N^{
-1} |\{1\leq i\leq N-K: (y_i,\ldots,y_{i+K}) \neq (x^N_i,\ldots,x^N_{i+K}) \}|
\\
\eps_N(\by) &=& N^{-1} (H_N(\by) - H_N(\bx^N))
\end{eqnarray*}
and then set
\begin{equation}
 \eps_N(\delta) =  \min \{\eps_N(\by) : \ \delta_N(\by) \geq \delta\} . \label{NK-eps}
 \end{equation}  
 We expect existence of a deterministic limit
\[ \eps(\delta) = \mbox{ a.s.-}\lim_{N \to \infty} \eps_N(\delta) . \]

\paragraph{A heuristic analysis} 
The purpose of this section is to give a heuristic argument for $\eps(\delta) \asymp \delta^2$.
Given $i$ and $l \geq K+1$,
consider the set
of sequences $\by$ such that
\begin{eqnarray*}
(y_j,\ldots,y_{j+K}) &=& (x^N_j,\ldots,x^N_{j+K}) \ \forall j \not\in [i+1,i+l] \\
(y_j,\ldots,y_{j+K}) &\neq & (x^N_j,\ldots,x^N_{j+K}) \ \forall j \in [i+1,i+l] .
\end{eqnarray*}
Over this set, let $D_{i,l}$ be the minimum of
$H_N(\by) - H_N(\bx^N)$
and let $\by^{(i,l)}$ be the minimizing sequence.
The distribution of $D_{i,l}$ essentially depends only on $l$,
not on $i$ or $N$; write $f_l(0+)$ for its density at $0+$.
Let's assume
\begin{equation}
\sum_{l \geq K+1} l^2 f_l(0+) = A < \infty .
\label{NK-ansatz}
\end{equation}

It is intuitively clear how to choose a sequence $\by$ which
minimizes $\eps_N(\by)$
for a given $\delta$.
Just fix a small $\eta > 0$,
and create a sequence of ``excursions" away from $\bx^N$ as follows.
For each pair $(i,l)$ such that $D_{i,l} < \eta l$,
choose $\by$ to equal $\by^{(i,l)}$ on %\footnote{overlaps will be negligible}
the sites $[i+K+1,i+l]$; set $\by = \bx^N$ elsewhere.   
See Figure 1.

\vspace{0.2in}

\centerline{
$\begin{array}{llllll}
01100011&\hspace*{-0.14in}010&\hspace*{-0.14in}0101110100010011010&\hspace*{-0.14in}10101101111&\hspace*{-0.14in}000101011010&\bx^N\\
01100011&\hspace*{-0.14in}\underline{011}&\hspace*{-0.14in}0101110100010011010&\hspace*{-0.14in}\underline{10001110100}&\hspace*{-0.14in}000101011010&\by
\end{array}$
}

\vspace{0.08in}
\centerline{
{\bf Figure 1.}
{\small 
Excursions of lengths $l = 3$ and $11$.  
Here $K = 2$.}}

\vspace{0.13in}
\noindent

With this scheme, $\delta$ will be the mean length of
possible excursions starting from a given site, that is
\[ \delta \sim \sum_{l \geq K+1} l \cdot \eta l f_l(0+) . \]
And $\eps$ is the mean increment of $H_N$ associated with possible excursions
starting from a given site, that is
\[
\eps \sim \sum_{l \geq K+1} (\eta l/2) \cdot \eta l f_l(0+) . \]
In other words
$\delta \sim A \eta, \ \eps \sim A \eta^2/2$,
giving
$\eps \sim (2A)^{-1} \delta^2$
which is the desired ``scaling exponent $=2$" result.

Why should the assumption (\ref{NK-ansatz}) be true?
Well, for large $l$ we expect central limit
behavior:
$D_l \approx {\rm Normal}(\mu l, \sigma^2 l)$
for some $\mu > 0$ and $0<\sigma^2 < \infty$.
This in turn suggests that $f_l(0+)$ should decrease at least geometrically
fast in $l$.

Note that the optimizing $\by^N$ in (\ref{NK-eps}) will have 
(in the $N \to \infty$ limit) some distribution $L_\delta$ of excursion lengths. 
The heuristic argument predicts that as $\delta \downarrow 0$ we have $L_\delta \cd L$ 
where the limit distribution has 
 $\Pb (L = l) \propto l f_l(0+)$ and 
$\Eb L < \infty$.

Simulations (Table 1) with $K=3$ are consistent with both the predicted scaling exponent $2$
and the prediction of existence of a $\delta \downarrow 0$ limit distribution $L$ for excursion lengths.  
Making a rigorous proof seems difficult and so we turn to a simpler example.

\vspace{0.3in}

\centerline{
$\begin{array}{clccc}
\theta&\ \ \ \delta&\eps&\eps/\delta^2&\Eb L_\delta\\
0.002&0.0397&4.85 \cdot 10^{-5}&0.0308&10.9\\
0.004&0.0774&2.00 \cdot 10^{-4}&0.0334&11.0\\
0.008&0.147&7.69 \cdot 10^{-4}&0.0354&11.3\\
0.016&0.266&2.75 \cdot 10^{-3}&0.0388&11.8
\end{array}
$}

\vspace{0.07in}
{\bf Table 1.}
Monte Carlo simulations with $K=3, N = 10,000$; 1000 repeats.
These are exact optimizations done by introducing a Lagrange
multiplier $\theta$ which penalizes matching $(K+1)$-tuples.
We find $c_3 = 0.3065$.

\subsection{Main model and results}
\label{sec-4}
Let $(\xi_i, i \geq 1)$ be i.i.d. copies of a strictly positive random variable $\xi$,
and write $G(x) = \Pb (\xi \leq x)$.  
Define the {\em benefit} function
\begin{equation}
f_n(A) = \left( |A| - \sum_{i=1}^{n-1} \xi_i \ind(i \in A, i+1 \in A) \right), 
\quad A \subseteq \{1,2,\ldots,n\} 
\label{Mn-def}
 \end{equation}
 where $\ind(B) = \ind_B$ denotes the indicator random variable associated with an event $B$.   
 Intuitively, we choose a set $A$ of items, getting reward $1$ from each item chosen but paying cost $\xi_i$ if we choose both $i$ and $i+1$; we seek to maximize 
 benefit = reward - cost.  
So we shall study
\begin{equation}
\label{eq:Mn}
M_n:= \max_{A \subseteq \{1,2,\ldots,n\}} 
f_n(A) .
 \end{equation}
To simplify exposition we will assume 
 \begin{equation}
 \mbox{$G$ has bounded continuous density $g$ with $g(\sfrac{1}{2}) > 0$}
 \label{G-assume}
 \end{equation}
 which implies
  \begin{equation}
0 < G(\sfrac{1}{2}) < 1   \label{G-assume-2}
 \end{equation}
 though we suspect that Theorems \ref{th1} and \ref{th2} remain true under some much weaker non-degeneracy assumptions.  See section \ref{sec-G-assume} for further remarks.

We will first prove the following:
\begin{Theorem}
\label{th1}
There exists $\sfrac{1}{2} \leq c \leq 1$, such that almost surely and in $L^1$, 
$$
\lim_{n \to \infty} n^{-1} M_n  = c .
$$
The constant $c$ is given by the forthcoming formula (\ref{eq:c}). If $\xi$ is an exponential random variable with parameter $\lambda>0$ then
$$
c = (1 - e^{-\lambda})^{-1} - \lambda^{-1}.
$$
\end{Theorem} 
We record the explicit value of $c$ only in the exponential case,
but one could use formula (\ref{eq:c}) to obtain explicit values for
other standard distributions.

We now formalize the setup in the introduction. 
The optimization problem (\ref{eq:Mn})  has a solution, a random subset $\Aopt_n \subseteq \{1,2,\ldots,n\}$, and Corollary \ref{C-unique} will show the solution is unique with probability $\to 1$ as $n \to \infty$.     Define the random variable:
\begin{equation}
\eps_n(\delta) :=
\min \left\{
n^{-1} (f_n(\Aopt_n) - f_n(B)): |B \bigtriangleup\Aopt_n| \geq
\delta n \right\}
 \label{eps-def}
\end{equation}
where the minimum is over all subsets $B \subset \{1,\cdots,n\}$ such that the symmetric difference with $\Aopt_n$ is at least $\delta  n $.  Our main result is the following:
\begin{Theorem}
\label{th2}
$\bar{\eps}(\delta) : = \lim_{n} \Eb \eps_n(\delta) $ exists
for all $0<\delta<1$, and 
\begin{equation}
\label{upper-epsilon-1}
\limsup_{\delta \downarrow 0} \delta^{-2}
\bar{\eps}(\delta) < \infty, \end{equation}
\begin{equation}
\label{lower-epsilon} \liminf_{\delta \downarrow 0} \delta^{-2} \bar{\eps}(\delta) > 0 . \end{equation}
\end{Theorem}
We now outline the key ideas in the proof, and the organization of the paper.
\begin{itemize}

\item Dynamic programming over i.i.d. data is essentially just study of a related Markov chain 
(section \ref{subsec:dp}),
and in our model there are simple 
{\em inclusion criteria}
for whether item $i$ is in the optimal solution.
The inclusion criterion involves  two Markov chains 
(one looking left, one looking right) and the cost $\xi_i$ 
(Table 2 and % Proposition \ref{z}).
Lemma \ref{LXn}).
\item By considering the related infinite-time 
{\em stationary} Markov chain and using the same inclusion criteria,
we can define a random subset 
$\Aopt \subset \Zbold$ 
interpretable as the solution of an infinite optimization problem
(section \ref{sec-stat}).
\item The $n \to \infty$ limit benefit in Theorem \ref{th1}
is just the mean benefit per item using $\Aopt$ in the infinite problem 
(section \ref{sec-limitc}).
\item Study of $\eps_n(\delta)$ is an 
``optimization under constraint" 
problem, most naturally handled via introduction of a Lagrange multiplier $\theta$.
So the $\Bnopt$ attaining the maximum in (\ref{eps-def}) 
can be studied as above by introducing a more complicated Markov chain 
parametrized by $\theta$ 
(section \ref{sec-quint}),
finding the inclusion criteria (Table 3),
formulating the parallel optimization under constraint problem, 
and observing that 
$\bar{\eps}(\delta)$ is representable via functions 
$\delta(\theta), \eps(\theta)$
defined in terms of the stationary distribution of the 
more complicated Markov chain 
(Proposition \ref{le:nMn}).
\end{itemize}
Without trying to write details, it seems intuitively clear that the 
methodology above could be implemented in more general dynamic programming models 
such as the NK model of section \ref{sec-NK}.
However, to complete the argument we need to analyze the $\theta \to 0$ 
behavior of the functions $\delta(\theta), \eps(\theta)$.
Even in our simple model, 
we do not have any useful explicit expression for the needed stationary distribution,
so we proceed via inequalities rather than using the exact formulas.
For the upper bound (\ref{upper-epsilon-1}) 
we just identify a ``local configuration" which can be replaced by a different local configuration at small extra cost
(section \ref{sec:upperbound}).
For the lower bound, we decompose the process into blocks by breaking  at certain special configurations, and then get bounds on the chance  that $\Bnopt$ differs from $\Anopt$ on a block 
%(section \ref{subsec:contd}) 
and bounds on the mean decrease in benefit if it does differ (section \ref{sec:mu}).  But these arguments rely on the particular combinatorial structure of our special model.  
It is not clear how readily they can be extended to general models.

\section{Analysis of optimal solutions}
\label{sec:triple}

\subsection{Non-uniqueness}
In the case $n=2$, if $\xi_1 > 1$ then  both $\{1\}$ and $\{2\}$ attain the maximum value $1$ 
of the optimization problem (\ref{eq:Mn}): the  optimizing set is not unique.
Corollary \ref{C-unique} shows that, provided some $\xi_i + \xi_{i+1} < 1$ is less than $1$, 
the optimizing set $\Anopt$ is unique, and by assumption (\ref{G-assume-2})  this proviso holds with probability $\to 1$ as $n \to \infty$.  
After this section we generally ignore the possibility of non-uniqueness.

We start with some terminology that will also be used later.  
For an integer interval $[g,d]$ with $d-g+1$ even, the two 
{\em complementary alternating subsets} $A_1, A_2$ are as shown in Figure 2.

\vspace{0.1in}

\centerline{
$
\begin{array}{cccccccccccc}
&g&-&(g+1)&-& &&&-&(d-1)&-&d\\
A_1&\circ&-&\bullet&-&\circ&\dots&\bullet&-&\circ&-&\bullet\\
A_2&\bullet&-&\circ&-&\bullet&\dots&\circ&-&\bullet&-&\circ\\
\end{array}
$}

\vspace{0.1in}

\centerline{
{\bf Figure 2.}  Included items marked $\bullet$, excluded items marked $\circ$.
}

\begin{Lemma}
\label{L-unique}
Let $n \geq 2$.  For almost all realizations of $\xi_1,\ldots,\xi_{n-1}$, the following are equivalent.\\
(a) The subset maximizing (\ref{eq:Mn}) is not unique. \\
(b) $n$ is even and the only optimal solutions are the two complementary alternating subsets of $[1,n]$.\\
(c) $n$ is even and $M_n = n/2$.
\end{Lemma}
{\bf Proof.} Either of (b,c) implies (a), so it is enough to show (a) implies (b) and (c).
Suppose distinct subsets $B_1$ and $B_2$ attain the maximum.  
Then a.s. the values of $\xi_i$ used in the optimal sum are identical, that is 
\begin{equation}
\{i: \ (i,i+1) \subset B_1\} =  \{i: \ (i,i+1) \subset B_2\} := \SS, \mbox{ say.} 
\label{BBS}
\end{equation}  
First suppose $\SS$ is empty.  
Then each of $B_1$ and $B_2$ has only isolated elements.
But amongst such sets, the maximum of (\ref{eq:Mn}) is attained 
(for $n$ odd) uniquely by the alternating subset giving $M_n = (n+1)/2$, or 
(for $n$ even) only by the complementary alternating subsets.
So $\SS$ empty implies (b) and (c).  
For  general $\SS$, take some $i \in B_1 \bigtriangleup B_2$, and then take the maximal interval $i \in [g,d] \subset [1,n]$ which is disjoint from $\SS$.  
Repeating the argument above, the restrictions of $B_1$ and $B_2$ to $[g,d]$ must be  
complementary alternating subsets.  
If $[g,d] \neq [1,n]$ then either $d+1$ or $g-1$ is in $\SS$ -- say $d+1$ -- and so $d+1 \in B_1\cap B_2$.
But exactly one of $B_1,B_2$ contains $d$, contradicting (\ref{BBS}).  
So $[g,d] = [1,n]$ and so $\SS$ is empty.  \qed

\begin{Corollary}
\label{C-unique}
If $\xi_i + \xi_{i+1} < 1$ for some $1 \leq i \leq n-2$ then a.s. $\Anopt$ is unique.
\end{Corollary}
{\bf Proof.} Fix $i$ with $\xi_i + \xi_{i+1} < 1$ and let $B$ be the alternating subset of $[1,n]$ containing $i$ and $i+2$.  Replacing $B$ by $B \cup \{i+1\}$ increases the benefit by $1 - \xi_i - \xi_{i+1} > 0$, 
so $B$ cannot be optimal, and the result follows from Lemma \ref{L-unique}.
\qed

\subsection{Dynamic programming}
\label{subsec:dp}
 
 Finding the maximum value and the maximizing subset of (\ref{eq:Mn}) is algorithmically easy by dynamic programming, as follows.
 Define
 \begin{eqnarray}
 V^L_{n,i} &=& \max_{i \in A \subseteq \{1,\ldots,i-1,i\}} 
\left( |A| - \sum_{j=1}^{i-1} \xi_j \ind(j \in A, j+1 \in A) \right)
\label{VLn}\\
 W^L_{n,i} &=& \max_{i \not\in A \subseteq \{1,\ldots,i-1,i\}} 
\left( |A| - \sum_{j=1}^{i-1} \xi_j \ind(j \in A, j+1 \in A) \right)
\label{WLn}
\end{eqnarray}
which differ in that the former requires $i \in A$ and the latter requires $i \not\in A$. 
The superscripts $L$ here and $R$ later indicate {\em left} and {\em right}.  
Note that in fact  $V^L_{n,i},  W^L_{n,i}$ above and  $ X^L_{n,i}$ below do not depend on $n$, but the notation is useful to distinguish from the limit process $X^L_i$ later.

From (\ref{VLn},\ref{WLn}) we see
 $V^L_{n,1} = 1, W^L_{n,1} = 0$ and by induction over 
$1 \leq i$
\begin{eqnarray*}
V^L_{n,i+1} &=& 1 + \max(V^L_{n,i} - \xi_{i},W^L_{n,i})
\\
W^L_{n,i+1} &=&  \max(V^L_{n,i}, W^L_{n,i}) 
\end{eqnarray*}
the two terms in the {\em max} indicating the choice of using or not using element $i$.  
Then $M_n = \max(V^L_{n,n},W^L_{n,n})$ and by examining which {\em max}
term is used at each stage leading to $M_n$ we can recover the
optimizing subset $\Anopt$.

We now describe an alternative, more useful way to obtain $\Anopt$. First, consider the evolution rule for the process
\begin{equation}
 X^L_{n,i} := V^L_{n,i} - W^L_{n,i} \label{XRn}
 \end{equation}
as $i$ increases; the rule is
\begin{eqnarray}
X^L_{n,i+1} &=& 1 + \max(0,X^L_{n,i} - \xi_{i}) - \max(0, X^L_{n,i}) \nonumber \\
&=& 1 + \max(-X^L_{n,i},-\xi_{i})\ind(X^L_{n,i} \geq 0) . \label{XL-7}
\end{eqnarray}
One can check by induction that $0 \leq X^L_{n,i} \leq 1$ and thus rewrite the recursion as
\[
X^L_{n,i+1} = \max(1-X^L_{n,i}, 1-\xi_{i}) 
. \]
For $n$ fixed we define the right processes analogously
 \begin{eqnarray}
 V^R_{n,i} &=& \max_{i \in A \subseteq \{i,i+1,\ldots,n\}} 
\left( |A| - \sum_{j=i}^{n-1} \xi_j \ind(j \in A, j+1 \in A) \right)
\label{VRn}\\
 W^R_{n,i} &=& \max_{i \not\in A \subseteq \{i,i+1,\ldots,n\}} 
\left( |A| - \sum_{j=i}^{n-1} \xi_j \ind(j \in A, j+1 \in A) \right),
\label{WRn}
\end{eqnarray}
with $V^R_{n,n}=1, W^R_{n,n}=0$.  Observe that the evolution rule for
the process
\begin{eqnarray}
X^R_{n,i} :=V^R_{n,i}- W^R_{n,i}  
\label{XRVW}
\end{eqnarray}
as $i$ decreases does not depend on $n$. In fact, we have
\begin{equation}
X^R_{n,i-1} = \max(1-X^R_{n,i}, 1-\xi_{i-1}). 
\label{XRni}
\end{equation}
The point is that we can determine the optimizing random set $\Anopt$
in terms of the quantities above.
Fix $i$ and consider the quantities 
$(X^L_{n,i},V^L_{n,i},W^L_{n,i})$, $\xi_i$, $(X^R_{n,i+1},V^R_{n,i+1},W^R_{n,i+1})$ 
and drop subscripts.  
We have four choices of whether to include (marked as $\bullet$ in Table 2)
or exclude (marked as $\circ$ in Table 2) items $i$ and $i+1$ in the optimal set
$\Anopt$.
For each choice, the table shows the absolute benefit of that choice, then the relative 
benefit (relative to the choice to exclude both items).
For each $i$ the optimal $\Anopt$ will contain, in positions $(i,i+1)$, the combination with largest 
relative benefit, and the final column indicates the criteria for use of each combination.

\vspace{0.2in}

\centerline{$\begin{array}{lccc}
- i - (i+1) - &
\mbox{absolute benefit}&
\mbox{relative benefit}&
\mbox{when used}\\
&&&\\
- \bullet - - \bullet - &
V^L + V^R - \xi &
X^L + X^R - \xi &
\mbox{ if } \xi < \min(X^L,X^R)\\
- \bullet - - \circ - &
V^L + W^R &
X^L &
\mbox{ if } X^R < \min(X^L,\xi)\\
- \circ - - \bullet - &
W^L + V^R &
X^R &
\mbox{ if } X^L < \min(X^R,\xi)\\
- \circ - - \circ - &
W^L + W^R &
0 &
\mbox{ never. } 
\end{array}$
}

\vspace{0.1in}
\centerline{{\bf Table 2.}  Inclusion criteria for $i, i+1$ in $\Anopt$.
}

\vspace{0.2in}
\noindent
(The case of non-uniqueness of $\Anopt$, Lemma \ref{L-unique}, is the case where 
$X^L_i$ and $X^R_i$ alternate between $0$ and $1$ throughout the interval $[1,n]$, and where we have equalities $X^L_i = X^R_{i+1} < \xi_i$.  Outside this case, one of the three strict inequalities holds.  We ignore the non-uniqueness possibility in the summary below.)

We summarize the argument above as follows.
\begin{Lemma}
\label{LXn}
For each $n$ define 
$X^L_{n,i}, 1 \leq i \leq n$
and
$X^R_{n,i}, 1 \leq i \leq n$
by
\begin{eqnarray}
X^L_{n,1} = 1; &&\quad X^L_{n,i+1} = \max(1-X^L_{n,i}, 1-\xi_{i}) , 1 \leq i \leq n-1 \label{XLn-rec}\\
X^R_{n,n} = 1; &&\quad X^R_{n,i-1} = \max(1-X^R_{n,i}, 1-\xi_{i-1}), 2 \leq i \leq n . \label{XRn-rec}
\end{eqnarray}
Then $\Anopt$ is the random subset of $\{1,2,\ldots,n\}$ specified by: 
for each $1 \leq i \leq n-1$,
\begin{eqnarray*}
\mbox{ if } \xi_i < \min(X^L_{n,i},X^R_{n,i+1}) 
& \mbox{then} & i \in \Anopt, \ i+1 \in \Anopt \\
\mbox{ if } X^R_{n,i+1} < \min(X^L_{n,i},\xi_i)
& \mbox{then} & i \in \Anopt, \ i+1 \not\in \Anopt \\
\mbox{ if } X^L_{n,i} < \min(X^R_{n,i+1},\xi_i)
& \mbox{then} & i \not\in \Anopt, \ i+1 \in \Anopt .
\end{eqnarray*}
\end{Lemma}
Let us emphasize two points:
\begin{itemize}
\item whether or not $i \in \Anopt$ depends only on the three r.v.s $X^L_{n,i}, \xi_i, X^R_{n,i+1}$
\item the only place where the value of $n$ enters is as the boundary condition $X^R_{n,n} = 1$.
\end{itemize}
In the next section, we show how
 to define a corresponding stationary process $((X^L_i,\xi_i,X^R_{i+1}), \ -
\infty < i < \infty)$.  By
applying the specification in Lemma \ref{LXn}  to this process, we will  define a set $\Aopt\subseteq \Z$ which will be shown (Lemma \ref{Lwc}) to be the limit of $\Anopt$.  As a consequence, we will be able to derive the limit of $M_n/n$.

\subsection{A stationary Markov chain and the infinite limit problem}
\label{sec-stat}
The recursion (\ref{XLn-rec}) specifies a Markov chain on the continuous state space $[0,1]$ 
with transitions
\begin{equation}
 x \to \max(1-x, 1 - \xi) , \label{kernel-1}
 \end{equation}
 where $\xi$ has distribution function $G$.
Write 
$F(x) = \Pb (X^L \leq x)$ for a stationary distribution function for this chain.
Then
\begin{eqnarray*}
F(x) &=& \Pb (\max(1-X^L,1-\xi) \leq x)\\
&=&\Pb (\min(X^L,\xi) > 1-x)\\
&=&\oG(1-x)\oF(1-x) 
\end{eqnarray*} 
where for any distribution function $F$ we write $\oF(x) = 1 - F(x)$.
Iterating this identity once gives
\[ F(x) = \oG(1-x) \left(1 - \oG(x)\oF(x) \right) \]
and solving this equation gives
\begin{equation}
\label{dist-1}F(x) = \frac{\oG(1-x)\oG(x)}{1 - \oG(x)\oG(1-x)}.
\label{dens-1}
\end{equation}
The assumption (\ref{G-assume}) that $G$ has a density implies that $F$ has a density, so in what follows we do not need to distinguish carefully between weak and strict inequalities for random variables with these distributions.

Now consider the infinite line graph, with vertices 
$- \infty < i < \infty$ and with i.i.d. edge-costs $\xi_i$ on edge
$(i,i+1)$ such that $\Pb (\xi_0+\xi_1<1)>0$, which is ensured by the condition $\oG(1/2)<1$.

\begin{Lemma}
\label{L1}
The recursion
\begin{equation}
X^L_{i+1} = \max(1-X^L_{i},1-\xi_i), \ - \infty < i < \infty 
\label{xL}
\end{equation}
defines uniquely a joint distribution for 
$((\xi_i,X^L_i), - \infty < i < \infty)$ 
in which $(X^L_i)$ is the stationary Markov chain with transition kernel (\ref{kernel-1}) and
stationary distribution (\ref{dist-1}).  And
\begin{equation} X^L_{i} = \phi(\dots,\xi_{i-2},\xi_{i-1}) \label{phi-def}
\end{equation}
for a certain function $\phi$ not depending on $i$.
\end{Lemma}
{\bf Proof.}
Having proved existence and uniqueness of the stationary distribution at (\ref{dens-1}), it only remains to
prove the measurability property (\ref{phi-def}).
Iterating (\ref{xL}) once shows
\begin{equation}
1 - \xi_i \leq X^L_{i+1} \leq \max(1-\xi_i, \xi_{i-1}) . \label{24.5}
\end{equation}
So outside the event 
$\{1 - \xi_i < \xi_{i-1}\}$ 
the value of $X^L_{i+1}$ depends only on 
$(\xi_{i-1},\xi_{i})$ 
and not on the value of $X^L_{i}$.
%\marginpar{M changed\\ $X^L_{i-1}\to X^L_{i}$}
So inductively on $Q \geq 1$ there exists a measurable function 
$\phi_Q$ such that
\[ X^L_1 = \phi_Q(\xi_{-2Q-1},\xi_{-2Q},\ldots,\xi_{0}) 
\mbox{ outside } 
\cap_{q=-Q}^{0} 
\{1 - \xi_{2q} < \xi_{2q-1}\} .
\]
Now (\ref{phi-def}) follows because 
$\Pb \left(  \cap_{q=-Q}^{0} 
\{1 - \xi_{2q} < \xi_{2q-1}\} \right) = \left( \Pb (\xi_0 + \xi_1 > 1) \right)^{Q+1} \to 0$.
\qed

If we define an ``$i$ decreasing" process by
\begin{equation}
 X^R_{i} = \phi(\dots,\xi_{i+2},\xi_{i+1},\xi_{i}) 
 \label{XRi}
 \end{equation}
then $(X^R_i)$ satisfies the analogous recursion
\begin{equation}
X^R_{i} = \max(1-X^R_{i+1},1-\xi_{i}), \ - \infty < i < \infty 
\label{XR-rec}
\end{equation}
and is distributed as the same stationary Markov chain.
Hence we have a rigorous definition 
of a unique (in distribution) stationary process 
$((X^L_i,\xi_i,X^R_{i+1}), \ - \infty < i < \infty)$ 
satisfying (\ref{xL},\ref{XR-rec})
which we will call the {\em triple process}.  
Note that from (\ref{phi-def},\ref{XRi})
\begin{equation}
\mbox{for each $i$ the three r.v.'s $X^L_i,\xi_i,X^R_{i+1}$ are independent.}
\label{3-ind}
\end{equation}
\begin{Lemma}
\label{LXinfty}
Let  $(X^L_i,\xi_i,X^R_{i+1}), \ - \infty < i < \infty$ be the stationary triple process.
Then there is a random subset $\Aopt$ of $\Z$ specified by: 
for each $-\infty < i < \infty $,
\begin{eqnarray*}
\mbox{ if } \xi_i < \min(X^L_{i},X^R_{i+1}) 
& \mbox{then} & i \in \Aopt, \ i+1 \in \Aopt \\
\mbox{ if } X^R_{i+1} < \min(X^L_{i},\xi_i)
& \mbox{then} & i \in \Aopt, \ i+1 \not\in \Aopt \\
\mbox{ if } X^L_{i} < \min(X^R_{i+1},\xi_i)
& \mbox{then} & i \not\in \Aopt, \ i+1 \in \Aopt .
\end{eqnarray*}
\end{Lemma}
{\bf Proof.} 
We need only check that the definition of $\Aopt$ is consistent, in that the
criterion for item $i$ to be excluded should be the same whether we look at
the pair $(i,i+1)$ or the pair $(i-1,i)$. 
(Of course this is intuitively clear from the consistency in the finite setting of Lemma \ref{LXn}, but let us give an algebraic verification anyway.) 
We need to check
\[ 
\{X^L_i < \min(X^R_{i+1},\xi_i)\}
\stackrel{?}{=}
\{X^R_i < \min(X^L_{i-1},\xi_{i-1})\} . \]
Using the recursions (\ref{XR-rec},\ref{xL}) for $X^R_i$ and $X^L_i$, we need to check
\[ 
\{\max(1-X^L_{i-1},1-\xi_{i-1}) < \min(X^R_{i+1},\xi_i)\}
\stackrel{?}{=} 
\{\max(1-X^R_{i+1},1-\xi_i) < \min(X^L_{i-1},\xi_{i-1})\} . \] 
But these are equal by applying the transformation $u \to 1-u$ to the right side.
\qed

Because the rule defining $\Aopt$ is translation-invariant, the augmented triple process
\[
((X^L_i,\xi_i,X^R_{i+1},\ind(i \in \Aopt)), \ - \infty < i < \infty) 
\] 
is also stationary.  
The next lemma shows this process is the limit of the corresponding finite-$n$ process.
The mode of convergence  can be viewed as a very elementary case of 
{\em local weak convergence} \cite{me101} of random graphical structures.  
In words, it asserts that relative to a random time-origin the finite processes approximate the limit process.
\begin{Lemma}
\label{Lwc}
Let $U_n$ be uniform on $\{1,\ldots,n\}$.
As $n \to \infty$
\[ 
((X^L_{n,U_n+i},\xi_{U_n+i},X^R_{n,U_n+i+1}, \ind(U_n +i \in \Anopt)), \ - \infty < i < \infty) 
\] \[
\cd 
((X^L_i,\xi_i,X^R_{i+1},\ind(i \in \Aopt)), \ - \infty < i < \infty)
\]
where the left side is defined arbitrarily for 
$U_n+i \not\in \{1,\ldots,n\}$ 
and where convergence in distribution is with respect to the usual product topology on infinite sequence space.
\end{Lemma}
{\bf Proof.}
Because the $X$'s are bounded and the $\xi$'s are i.i.d., 
the sequence of processes is tight in the product topology. 
Write
\[
((\hat{X}^L_i,\hat{\xi}_i,\hat{X}^R_{i+1},\ind(i \in  \hat{A}^{\mbox{{\tiny opt}}} )), \ - \infty < i < \infty)
\]
for a subsequential weak limit.
Clearly $(\hat{\xi}_i) \ed (\xi_i)$.
Because for each $n$ the process 
$(X^L_{n,i},\xi_i)$ satisfies recursion (\ref{XLn-rec}),
the limit $(\hat{X}^L_i,\hat{\xi}_i)$ 
satisfies this recursion, and so by the
``uniqueness of joint distribution" assertion of Lemma \ref{L1},
$(\hat{X}^L_i,\hat{\xi}_i)  
\ed (X^L_i,\xi_i)$.
Applying the same argument to $X^R$ we deduce 
\[ 
((X^L_{U_n+i},\xi_{U_n+i},X^R_{n,U_n+i+1} ), \ - \infty < i < \infty) 
\cd 
((X^L_i,\xi_i,X^R_{i+1} ), \ - \infty < i < \infty) .
\]
For fixed $i_0$ the event  $i_0 \in \Aopt$ is a function of the limit process, the function implied by Lemma \ref{LXinfty}, and by a standard fact (\cite{MR0233396} Theorem 5.2) it is enough to check that this function is a.s. continuous with respect to the limit process.  
But this just requires that the probability of an {\em equality} between some two of  $X^L_{i_0},\xi_{i_0},X^R_{i_0+1}$ should be zero, which follows from their independence (\ref{3-ind}) and existence of densities (\ref{G-assume},\ref{dens-1}).
\qed

\iffalse
The result now follows because the rule for determining whether $i \in \Aopt$ 
is a.s. continuous.
More precisely, by the Skorokhod representation theorem, we can assume
that with probability one, we have
\begin{eqnarray}
\label{limas}((X^L_{U_n+i},\xi_{U_n+i},X^R_{n,U_n+i+1} ), \ - \infty < i < \infty) 
\to 
((X^L_i,\xi_i,X^R_{i+1} ), \ - \infty < i < \infty) .
\end{eqnarray}
We define $\tau^R(i) = \inf\{n\geq i, \xi_n+\xi_{n-1}\leq 1\}$ and
  $\tau^L(i)=\sup\{n<i, \xi_n+\xi_{n-1}\leq 1\}$. Then as shown in
  the proof of Lemma \ref{L1}, the $X^L_j,X^R_{j+1}$ for $\tau^L(i)
  \leq j \leq \tau^R(i)$ are functions of the $((\xi_j ), \tau^L(i)
  \leq j \leq \tau^R(i))$ only. Hence there exists $\Psi$ such that
\begin{eqnarray*}
\ind(i \in {\Aopt}) = \Psi(\xi_j , \tau^L(i)
  \leq j \leq \tau^R(i)).
\end{eqnarray*}
By (\ref{limas}), we have $((\xi_{U_n+j} ), \tau^L(i)
  \leq j \leq \tau^R(i))\to ((\xi_j ), \tau^L(i)
  \leq j \leq \tau^R(i))$ and hence for $n$ sufficiently large,
\begin{eqnarray*}
\ind(U_n+i \in {\Anopt}) = \Psi(\xi_{U_n+j} , \tau^L(i)
  \leq j \leq \tau^R(i)),
\end{eqnarray*}
which implies
\begin{eqnarray*}
\ind(U_n+i \in {\Anopt}) \to \ind(i \in {\Aopt}) \mbox{ a.s.}
\end{eqnarray*}
\fi

\subsection{Proof of Theorem \ref{th1}}
\label{sec-limitc}
Because
\[ M_n = \sum_{i=1}^n \ind(i \in \Anopt) - \sum_{i=1}^{n-1} \xi_i \ind(i \in \Anopt, i+1 \in \Anopt) \]
we can write
\[n^{-1}  \Eb M_n = \Pb(U_n \in \Anopt) - \Eb \xi_{U_n}\ind(U_n \in \Anopt, U_n+1 \in \Anopt) \ind (U_n \neq n) \] 
and then by Lemma \ref{Lwc}
\[ n^{-1}  \Eb M_n \to c:= \Pb (0 \in \Aopt) - \Eb \xi_0 \ind(0 \in \Aopt, 1 \in \Aopt) . \]
Note that clearly $c \leq 1$; the other inequality $c \geq 1/2$ holds because the subset 
$\{1,3,5,\ldots\}$ is a feasible choice.

We now exploit  the {\em method of bounded differences} \cite{macd91} in a very routine way. 
We observe that $M_n = m_n (\xi_1,\cdots,\xi_n)$ for a certain function $m_n$ with the property
\begin{quote}
changing any one argument of $m_n(z_1,\ldots,z_n)$ changes the value of $m_n(\cdot)$ by at most $1$
\end{quote}
This property holds because  $\Anopt$ will never contain a pair $(i,i+1)$ for which $\xi_i > 1$.
And this property implies the well-known Azuma-Hoeffding inequality of the form 
 (see e.g. \cite{tal96})
$$
\Pb ( | M_n - \med(M_n) | \geq   t ) \leq 4 \exp( - \sfrac{t^2}{4n}). 
$$
It is now routine to use this large deviation inequality to establish the a.s. and $L^1$ convergence 
of $n^{-1} M_n$ to $c$.

\iffalse
As a direct consequence of Lemma \ref{Lwc} we can determine the limit behavior of the optimal reward
$M_n$.  First note that by a standard concentration argument,  $n^{-1}M_n$ is close to its mean. Indeed, define the function $m_n$ from $\mathbb R_+^n$ to $\mathbb R$, by $m_n( z_1, \cdots, z_n) = \min_{A \subseteq \{1,\cdots,n\} } |A| - \sum_{i =1}^{n-1} z_i 1 ( i \in A, i+1 \in A)$. By definition, we have $M_n = m_n (\xi_1,\cdots,\xi_n)$. The function $m_n$ is clearly $1$-Lipschitz for the Hamming-distance: $d_H( (z_1, \cdots,z_n), (z'_1, \cdots,z'_n) ) = \sum_{i=1} ^ n 1 ( z'_i \ne z_i)$. Thus, from Talagrand's concentration inequality (see for example Proposition 2.1 in Talagrand \cite{new look at independence}):
$$
\Pb ( | M_n - \med(M_n) | \geq   t ) \leq 4 \exp( - \frac{t^2}{4n}). 
$$
Integrating over all $t$, we obtain: $ |\Eb M_n -  \med(M_n) | \leq \Eb | M_n - \med(M_n) | \leq 4 ( \pi n) ^{1/2}$. We deduce that for all $t \geq 4 \pi^{1/2} n^{-1/2}$, 
$$
\Pb  ( | M_n - \Eb M_n | \geq   t n  ) \leq 4 \exp( -  n  ( t - 4 \pi^{1/2} n^{-1/2} ) ^2 / 4). 
$$
It follows that $n^{-1} M_n  $ converges almost surely and in norm $L^1$ to $\lim_n  n^{-1}  \Eb M_n$ (provided it exists). 
\fi

To evaluate $c$, abbreviate $(X^L_0,\xi_0,X^R_1)$ to $(X^L,\xi,X^R)$ and use the Lemma \ref{LXinfty}  definition of $\Aopt$ to write
\begin{eqnarray*}
 \Pb (0 \in \Aopt) &=& 1 - \Pb (X^L < \min(X^R,\xi)) \\
&=& 1 - \sfrac{1}{2}(1 - \Pb (\xi < \min(X^L,X^R)))\mbox{ by symmetry } \\
&=& 
\sfrac{1}{2} + \sfrac{1}{2} \Pb (\xi < \min(X^L,X^R)) 
\end{eqnarray*}
and then
\begin{equation}
\label{eq:c}
c = 
\sfrac{1}{2} + \sfrac{1}{2} \Pb (\xi < \min(X^L,X^R)) 
- \Eb \xi \ind(\xi < \min(X^L,X^R)) .
\end{equation}
Recall that $X^L$, $\xi$ and $X^R$ are independent and that $X^L$ and $X^R$ have common distribution $F$ given in terms of $G$ by (\ref{dist-1}).  
So (\ref{eq:c}) constitutes a formula for $c$ in terms of the underlying distribution function $G$ of $\xi$.

We now evaluate $c$ in the special case where $\xi$ has the exponential($\lambda$) distribution:
\[ \oG(x) = e^{-\lambda x}, \ 0 < x < \infty , \]
so that, from formula  (\ref{dist-1}), we have:
\[ F(x) = \frac{e^{-\lambda (1-x)} (1-e^{-\lambda x})}{1 -
  e^{-\lambda}} = \frac{e^{\lambda x} - 1}{e^\lambda -1}. \]
We deduce 
 \begin{eqnarray*}
\Pb (\xi < \min(X^L,X^R))
&=& \int_0^1 \lambda e^{-\lambda u} P^2(X^L>u) \ du \\
&=& \frac{\lambda}{(e^\lambda -1)^2} 
\int_0^1 e^{-\lambda u} \left(e^\lambda - e^{\lambda u}\right)^2 \ du \\
&=& \frac{\lambda}{(e^\lambda -1)^2} 
\left(
e^{2 \lambda} \int_0^1 e^{-\lambda u} du 
- 2 e^\lambda 
+ \int_0^1 e^{\lambda u} du \right)\\
&=& \frac{\lambda}{(e^\lambda -1)^2} 
\left(
\frac{e^{2\lambda} (1-e^{-\lambda})}{\lambda} 
- 2 e^\lambda 
+ \frac{e^\lambda -1}{\lambda} \right)  \\
&=& \frac{e^{2\lambda} - 2\lambda e^\lambda - 1}{
(e^\lambda - 1)^2} .
\end{eqnarray*}
\begin{eqnarray*}
 \Eb \xi \ind(\xi < \min(X^L,X^R)) 
&=& \int_0^1 u \lambda e^{-\lambda u} P^2(X^L >u) \ du \\
&=& \frac{\lambda}{(e^\lambda -1)^2} 
\int_0^1 u e^{-\lambda u} \left(e^\lambda - e^{\lambda u}\right)^2 \ du \\
&=& \frac{\lambda}{(e^\lambda -1)^2} 
\left(
e^{2 \lambda} \int_0^1 u e^{-\lambda u} du 
-  e^\lambda 
+ \int_0^1 u e^{\lambda u} du \right)\\
&=& \frac{\lambda}{(e^\lambda -1)^2} 
\left(
\frac{e^{2\lambda}% (1-e^{-\lambda})}{\lambda} 
(1 - (1+\lambda)e^{-\lambda})}{\lambda^2}
-  e^\lambda 
+ \frac{1 + (\lambda - 1)e^\lambda }{\lambda^2} \right) \\
&=& \frac{1}{\lambda (e^\lambda -1)^2}
\left( e^{2\lambda} - (\lambda^2 +2)e^\lambda + 1 \right) .
\end{eqnarray*}
Combining,
\begin{eqnarray*}
c &=& \frac{1}{2} + 
\frac{\frac{\lambda}{2} (e^{2 \lambda} - 2 \lambda e^\lambda - 1)
- (e^{2 \lambda} - (\lambda^2 +2)e^\lambda + 1)}
{\lambda (e^\lambda - 1)^2}\\
&=& \frac{1}{2} + 
\frac{
(\frac{\lambda}{2}-1)e^{2 \lambda} + 2e^\lambda - \frac{\lambda}{2} - 1}
{\lambda (e^\lambda - 1)^2}\\
&=& \frac{1}{2} + 
\frac{
\frac{\lambda}{2}(e^{2 \lambda} -1) - (e^\lambda - 1)^2}
{\lambda (e^\lambda - 1)^2}\\
&=& \frac{1}{2} + \frac{e^\lambda +1}{2(e^\lambda -1)} - \frac{1}{\lambda} \\
&=& \frac{1}{1-e^{-\lambda}} - \frac{1}{\lambda} .
\end{eqnarray*}

\section{The upper bound in Theorem \ref{th2} }
 \label{sec:upperbound}

Local weak convergence (Lemma \ref{Lwc} above and Lemma \ref{Lwc2} below) 
 provides one sense in which the $n \to \infty$ limit of the solution $\Anopt$ of the size-$n$ optimization problem is $\Aopt$. 
 A logically different sense is provided by  coupling, as follows. 
 Part of the stationary triple process is the doubly-infinite i.i.d. sequence 
 $(\ldots,\xi_{-1},\xi_0,\xi_1,\xi_2,\ldots)$.  
 For each $n$ use these same r.v.'s 
 $\xi_1,\ldots,\xi_n$ to construct $\Anopt$.
 Because of boundary effects it is not always true that 
 $\Aopt \cap [1,n] = \Anopt$.
 But we expect the sets to coincide ``away from the boundary", and Lemma \ref{lem1}(b) below 
 provides one expression of this equality.  
 We call this technique {\em localization}.

 \subsection{Optimality properties of $\Aopt$}

Lemma \ref{LXinfty} gave a concise definition of $\Aopt$ but did not explicitly identify its optimality properties.  Lemma \ref{lem1} below will relate  $\Aopt$ to certain finite optima and thereby allow us to deduce some explicit properties.

The benefit function $f_n(A)$ and its maximum value $M_n$ defined  at (\ref{Mn-def},\ref{eq:Mn})  
refer to subsets of $[1,n]$, and it is convenient to make the corresponding definitions for an arbitrary 
interval $[\ell,m]$:
\begin{eqnarray}
f_{[\ell,m]}(A)  &:=&
 |A| - \sum_{i=\ell}^{m-1} \xi_i \ind(i \in A, i+1 \in A) , 
\quad A \subseteq \{\ell,\ell+1,\ldots,m\}
\label{Mlm}\\
M_{[\ell,m]} &:=& \max_{A \subseteq \{\ell,\ell+1,\ldots,m\}}  
f_{[\ell,m]}(A)  
\end{eqnarray}
and denote by $\Aopt_{[\ell,m]}$ the corresponding optimizing set.

\begin{Lemma}
\label{lem1}
(a)  If $\xi_{i-1}+\xi_ i \leq 1$ then $i \in \Aopt$. \\ 
(b) If $\ell < m$ and $\xi_{\ell-1}+\xi_ \ell \leq 1$ and $\xi_{m-1}+\xi_ m \leq 1$ then $\Aopt_{[\ell,m]}$ is unique and
\begin{eqnarray}
\label{eq:finite}\Aopt\cap [\ell,m]=\Aopt_{[\ell,m]}.
\end{eqnarray}
If furthermore $[\ell, m] \subseteq [1,n]$ then $\Anopt\cap [\ell,m]=\Aopt_{[\ell,m]}$ 
(interpreting $\xi_0 = 0$ if $\ell = 1$).  \\
(c) If both $i, i+1\in \Aopt$ then $\xi_i\leq 1$. \\
(d) If $\xi_i+\xi_{i+1}>1$ then $i,i+1$ and $i+2$ together cannot belong
to $\Aopt$. \\
(e) Let $k \geq 2$.  If $[g, g+2k - 1]$ is an interval such that
$\xi_g > \xi_{g+1} > \ldots > \xi_{g+2k-1} >  \xi_{g+2k}$ and 
\[  \xi_{j}+\xi_{j+1}>1, \quad g \leq j \leq g  +2k - 2 \]
 then  $\Aopt \cap [g,g+2k -1]$ must be one of the two 
 complementary alternating sequences in $[g, g+2k- 1]$.
\end{Lemma}
{\bf Proof.}
(a)
 If $\xi_{i-1}+\xi_ i \leq 1$ then
$X^L_i\geq 1-\xi_{i-1}\geq \xi_i$
and hence from the Lemma \ref{LXinfty} definition
we see that 
for any possible value of $X^R_{i+1}$, we have $i \in\Aopt$.

(b) 
First note that both $\ell$ and $m$ are in $\Aopt_{[\ell,m]}$, for otherwise adding each element would increase $f_{[\ell,m]}(\Aopt_{[\ell,m]})$ by at least  $1-\xi_\ell$ and $1-\xi_{m-1}$ respectively.
Next note that by rewriting Lemma \ref{LXn} 
(which concerns the special case $[\ell,m] = [1,n]$)
for general $[\ell,m]$, we have a construction of $\Aopt_{[\ell,m]}$
in terms of processes 
$X^L_{[\ell,m],i}$ and $X^R_{[\ell,m],i}$ for $\ell\leq i\leq m$ 
defined by the recursions analogous to (\ref{XLn-rec},\ref{XRn-rec}).
By (a) both $\ell$ and $m$ are in $\Aopt$.
We have now shown $X^L_{[\ell,m],\ell}=1-\xi_{\ell-1} =X^L_\ell$ and $X^R_{[\ell,m],m-1}=1-\xi_{m-1} =X^R_{m-1}$; 
because the restricted and unrestricted processes have the same boundary conditions and satisfy the same recursions over $[l,m]$ they must agree throughout the interval.  
Finally, because both endpoints $\ell$ and $m$ are in $\Aopt_{[\ell,m]}$ it cannot fit the 
``complementary alternating sequences" criteria for non-uniqueness (Lemma \ref{L-unique}).  
The same argument works for $\Anopt$.

One could prove (c,d,e) algebraically 
from the definition of $\Aopt$, but it is more intuitive to exploit the finite 
optimality criterion as follows.
From assumption (\ref{G-assume-2}) 
there are infinitely many $\ell$ with $\xi_{\ell -1}+   \xi_\ell \leq 1$, 
and so for any given $i$ there is a (random) long interval 
$[\ell,m]$ containing $i$ 
for which by (b) 
$\Aopt\cap [\ell,m]=\Aopt_{[\ell,m]}$.
In other words, the restriction of $\Aopt$ to $[\ell,m]$ 
is the solution of the finite optimization problem (\ref{Mlm}), and we can derive its properties
by considering the effect of local changes.
Now (c) and (d) follow from the observations \\
(for (c)): 
if $i,i+1$ are in $\Aopt$ then removing $i+1$ will give a relative benefit of at least $\xi_i-1$.\\
(for (d)): 
if $i,i+1,i+2$ are all in 
$\Aopt$ then removing $i+1$ will give a relative benefit of
$\xi_i+\xi_{i+1} - 1$.

For (e), consider $j \in [g, g+2k]$.
By (d) we cannot have $\{j,j+1,j+2\} \subset \Aopt$. 
If $j$ and $j+1$ but not $j+2$ are in $\Aopt$ then deleting $j+1$ while adding $j+2$ would increase the benefit by at least $\xi_j - \xi_{j+2} > 0$, which is impossible.  
It follows that we cannot have $\{j,j+1\} \subset \Aopt$.  
Thus $\Aopt \cap [g,g+2k -1]$ contains only isolated elements.   
It is now easy to check that one can change $\Aopt \cap [g,g+2k -1]$ into one of the alternating sequences on $ [g,g+2k -1]$ in such a way that the cardinality does not decrease, and the end items 
$g, g+2k-1$ change (if at all) only from included to excluded.  
Thus the change can only increase the benefit; appealing to the uniqueness property (b) in a larger interval establishes (e).\qed

\subsection{Proof of upper bound}
\label{sec-pr-ub}
In this section we prove the  bound
\begin{eqnarray}
\label{upper-epsilon-2}
\limsup_{\delta \downarrow 0 } \delta^{-2} \limsup_n \Eb \eps_n (\delta) < \infty
\end{eqnarray} 
via a simple
construction of near-optimal sets.  We first describe a particular configuration.
%We call the configuration with an alternating pattern
%$\bullet\circ\bullet\circ\dots$ a switch. 
Let $g,d\in \Zbold$ such that $d-g =
2k$ for some $k \geq 2$, and consider the sets $A$ and $B$ below 
\begin{eqnarray*}
\begin{array}{cccccccccccccccc}
&g&-&(g+1)&-& &&&&&&(g+2k-1)&-&g+2k\\
A&\bullet&-&\circ&-&\bullet&-&\circ&\dots&\bullet&-&\circ&-&\bullet\\
B&\bullet&-&\bullet&-&\circ&-&\bullet&\dots&\circ&-&\bullet&-&\bullet\\
\end{array}
\end{eqnarray*}
where  $|A\bigtriangleup B| = 2k -1$ and
the difference between the benefits of $A$ and $B$ is:
\begin{eqnarray}
f_{[g,g+2k]}(A)-f_{[g,g+2k]}(B) &=& (k+1) -\left(
(k+2)-\xi_g-\xi_{g+2k-1}\right)\nonumber \\
&=& \xi_g+\xi_{g+2k-1}-1. \label{fggk}
\end{eqnarray}
Now fix $k \geq 2$ and
$\alpha>0$ such that $\alpha k<1/2$. Consider the event $\Omega_g$ defined by: 
\begin{eqnarray*}
\xi_g > \xi_{g+1}>\xi_{g+2}>\dots
>\xi_{g+2k-2}> \xi_{g+2k-1}>\sfrac{1}{2}>\xi_{g+2k}; \\
\xi_{g-1} + \xi_g < 1; \ 
\xi_{g+2k-1} + \xi_{g+2k} <1 .
\end{eqnarray*}  
By assumption (\ref{G-assume}) this event has non-zero probability.  
If this event occurs, Lemma \ref{lem1}(a) shows that $\Aopt$ contains $g$ and $g+2k$,
and then Lemma \ref{lem1}(e) implies that $\Aopt \cap [g,g+2k]$ is the set $A$ above. 
By applying Lemma \ref{lem1}(b) to $[\ell,m] = [g,g+2k]$ we have the analog in the finite $n$ setting:
if $\Omega_g$ occurs for $[g,g+2k] \subset [1,n]$ then $\Anopt \cap [g,g+2k]$ is the set $A$ above.  
So if we change $\Anopt$ by replacing pattern $A$ by pattern $B$ on such an interval, then from (\ref{fggk}) the decrease in benefit equals 
$ \xi_g+\xi_{g+2k-1} - 1 > 0$. 
Now define 
\[ \Omega_g^{(\alpha)} = \Omega_g \cap   \{1 < \xi_g+\xi_{g+2k-1} < 1 + 2k\alpha\} \]
\[ q(\alpha) = \Pb (\Omega_g^{(\alpha)}) \]
\[ r(\alpha)=\Eb (\xi_{g}+\xi_{g+2k - 1} - 1) \ind(\Omega_g^{(\alpha)})  \] 
so that $r(\alpha)$ is the unconditional mean increase in benefit from the possible change, now performed only if event $\Omega_g^{(\alpha)}$ happens.
Using assumption (\ref{G-assume}) we see that 
$ (\xi_{g}+\xi_{g+2k - 1} )$ restricted to $\Omega_g^{(\alpha)} $ has a continuous density which is non-zero at $1$, which easily implies
 that for fixed $k$
 \begin{equation}
 \label{qr-lim}
 q(\alpha) \sim \bar{q}\alpha, \quad r(\alpha) \sim \bar{r}\alpha^2 \mbox{ as } \alpha \downarrow 0
 \end{equation}
 for constants $\bar{q}, \bar{r} \in (0,\infty)$.
 
 Given $n$ and the optimal set $\Anopt$, construct a near-optimal set $B^{(\alpha)}_n$ as follows. 
 Let $g_1 = 1$ and let 
 \[ [g_1, g_1 + 2k], \ [g_2, g_2 + 2k], \ [g_3, g_3 + 2k], \ldots, [g_{j_n}, g_{j_n}+2k] \]
 be the adjacent disjoint intervals in $[1,n]$ containing $2k+1$ integers.  
 For each such $g = g_j$, if event $\Omega_g^{(\alpha)} $ occurs, then on $[g, g+2k]$ replace pattern $A$ by pattern $B$.
 
Letting $n\to\infty$ and using the weak law of large numbers, we get
\begin{eqnarray*}
\sfrac{1}{n} |B^{(\alpha)}_n\bigtriangleup \Aopt_n| &\to& 2kq(\alpha)/(2k+1) \mbox{ in probability,}\\
\sfrac{1}{n} (f_n(\Aopt_n)-f_n(B^{(\alpha)}_n)) &\to& r(\alpha)/(2k+1)  \mbox{ in probability.}
\end{eqnarray*}
If $\frac{1}{n} |B^{(\alpha)}_n\bigtriangleup \Aopt_n| \leq kq(\alpha)/(2k+1)$ then redefine $B^{(\alpha)}_n$ to be the empty set.   Then (taking $k =3$ for concreteness)
\begin{eqnarray*}
\sfrac{1}{n} |B^{(\alpha)}_n\bigtriangleup \Aopt_n| &\geq & 3q(\alpha)/7 \\
\lim_n \sfrac{1}{n} (\Eb f_n(\Aopt_n)- \Eb f_n(B^{(\alpha)}_n)) &=& r(\alpha)/7  .
\end{eqnarray*}
The upper bound (\ref{upper-epsilon-2}) now follows from the $\alpha \to 0$ asymptotics (\ref{qr-lim}).

\section{Proof of Theorem \ref{th2} : the lower bound}
 \label{sec:lowerbound}

\subsection{Analysis of near-optimal solutions: the quintuple process}
\label{sec-quint}
Throughout section \ref{sec:lowerbound}  we fix a constant $\tau> 0$ such that
\begin{equation}
\label{eq:deftau}
G \left( \sfrac{1}{2} - \tau \right) > 0.
\end{equation}
Such a constant exists  by assumption (\ref{G-assume}).  
To study near-optimal solutions, fix a Lagrange multiplier $\theta$ such that
\begin{equation}
0 < \theta < \tau . \label{theta-tau}
\end{equation}
%\marginpar{M added bar}
We will derive the existence of, and derive an exact expression for, the function $\bar{\eps}(\delta) = \lim_n \Eb \eps_n (\delta)$ when $\delta$ is sufficiently small.
  The expression is an implicit function representation $\bar{\eps}(\delta(\theta))=\eps(\theta)$ via two functions $\eps(\theta), \delta(\theta)$ defined (\ref{del-def-0},\ref{eps-def-0}) in terms of the stationary distribution of a certain {\em quintuple process}.

We study the modified optimization problem in which we get an extra reward $\theta$ for choosing an item which is not in $\Aopt_n$
or for not choosing an item which is in $\Aopt_n$:
\begin{equation}
 \max_{A \subseteq [n]}
\left( |A| - \sum_{i=1 }^{n} \xi_i \ind(i \in A, i+1 \in A) + \theta |A \bigtriangleup \Aopt_n|  \right)
 . \label{iop2}
 \end{equation}
To study this we modify (\ref{VLn},\ref{WLn}) to
 \begin{eqnarray}
 \bV^L_{n,i} = \max_{i \in A \subseteq \{1,2,\ldots,i\}} 
\left( |A| - \sum_{j=1}^{i-1} \xi_j \ind(j,j+1 \in A) 
+ \theta \left| (A \bigtriangleup \Aopt_n)\cap \{1,2,\ldots,i \} \right| \right) &&
\label{VR2}
\\
 \bW^L_{n,i} = \max_{i \not\in A \subseteq \{1,2,\ldots,i\}} 
\left( |A| - \sum_{j=1}^{i-1} \xi_j \ind(j,j+1 \in A) 
+ \theta \left| (A \bigtriangleup \Aopt_n)\cap \{1,2,\ldots,i \} \right| \right) && .
\label{WR2}
\end{eqnarray}
We also define $\bM_n = \max(\bV^L_{n,n}, \bW^L_{n,n})$ and write $\Bnopt$ for the corresponding optimizing set.  Note that these quantities depend on $\theta$. 
Analogous to the definition (\ref {XRn}) of $X^L_{n,i}$ we define
\begin{eqnarray*}
Z^L_{n,i} := \bV^L_{n,i}-\bW^L_{n,i}.
\end{eqnarray*}
Then as the analog of (\ref{XL-7}) we can obtain the recursion
\[Z^L_{n,i+1} =  1-\min(Z^L_{n,i}, \xi_i)\ind(Z^L_{n,i}>0)+\theta J_{n,i+1}\]
where
\begin{eqnarray*}
Z^L_{n,1} &=& 1+\theta J_{n,1} \\
J_{n,i} &=& \ind(i\notin \Aopt_n) - \ind(i\in \Aopt_n) .
\end{eqnarray*}

Recall from section \ref{sec-stat} the stationary triple process
$((X^L_i,\xi_i,X^R_{i+1}), - \infty < i < \infty)$ and define
 \[ J_i =\ind(i\notin \Aopt) - \ind(i\in \Aopt) . \] 
 Just as the stationary triple process is interpretable (Lemma \ref{Lwc}) as an $n \to \infty$ limit of the process 
 $(X^L_{n,i},\xi_i,X^R_{n,i+1})$, we want to define a process which will be the limit of 
 $(Z^L_{n,i},X^L_{n,i},\xi_i,X^R_{n,i+1})$.  
 So define a {\em quadruple process} $(Z^L_i,X^L_i,\xi_i,X^R_{i+1})$ to be a process such that\\
 (i) $(X^L_{i},\xi_i,X^R_{i+1})$ evolves as the triple process\\
 (ii) $Z^L_{i}$ satisfies the recursion
\begin{equation}
\label{eq:recZ}
Z^L_{i+1} = 1-\min(Z^L_{i}, \xi_i)\ind(Z^L_{i}>0)+\theta J_{i+1} . %\ - \infty < i < \infty 
\end{equation}
%and introduce the quadruple process:  $((Z^L_i,X^L_i,\xi_i,X^R_{i+1}), - \infty < i < \infty)$.  
Recall $0<\theta<\tau$.
\begin{Lemma}
\label{L4}
The quadruple process  $((Z^L_i,X^L_i,\xi_i,X^R_{i+1}), - \infty < i < \infty)$ has a unique stationary distribution, for which
\begin{equation} Z^L_{i} = \psi(\dots,\xi_{i-2},\xi_{i-1},\xi_i,X^R_{i+1}) \label{psi-def}
\end{equation}
for a certain function $\psi$ not depending on $i$.   On the event
$\{\xi_{i-1}+\xi_i\leq 1-\tau\}$, we have
\begin{eqnarray}
X^L_{i+1}=1-\xi_i, && Z^L_{i+1}=1-\xi_i+\theta J_{i+1}. \label{XxX}
\end{eqnarray}
\end{Lemma}
{\bf Proof.}  
Recursion  (\ref{eq:recZ}) implies 
$ Z^L_{i+1} \geq 1 - \xi_{i} + \theta J_{i+1} $. 
Thus iterating once (\ref{eq:recZ}) and using this last inequality, we obtain:
$$
1 - \xi_{i} + \theta J_{i+1} \leq Z^L_{i+1} \leq 1 - \min (1 - \xi_{i-1} + \theta J_{i}, \xi_i ) \ind(1 - \xi_{i-1} + \theta J_{i}  >0) + \theta J_{i+1} .
$$
Thus, on the event $\{\xi_{i-1} + \xi_{i} \leq 1 - \theta \}$ we have 
$
Z^L_{i+1} = 1 - \xi_i + \theta J_{i+1}
$ 
and also, by (\ref{24.5}), we have  $X^L_{i+1} = 1 - \xi_i$, establishing (\ref{XxX}). 
Assumption (\ref{G-assume}) implies that the event $\{\xi_{i-1}+\xi_i\leq 1-\tau\}$ occurs for infinitely many $i < 0$, so in particular 
$K:= \max\{i < 0: \xi_{i-1} + \xi_{i} \leq 1 - \tau \}$ is finite.  
By the recursion (\ref{eq:recZ}) we can write $Z^L_0$ in the form
\[ Z^L_0 = \psi^1 (\xi_{K+1}, \xi_{K+2},\ldots,\xi_{-1}; Z^L_{K+1}; J_{K+2}, J_{K+3},\ldots, J_0) \]
for some function $\psi^1$.  
Then by (\ref{XxX}) with $Z^L_i = Z^L_{K+1}$ we can rewrite as 
\[ Z^L_0 = \psi^2 (\xi_K,\xi_{K+1}, \xi_{K+2},\ldots,\xi_{-1}; J_{K+1}, J_{K+2}, J_{K+3},\ldots, J_0). \]
By the definition of $\Aopt$, each $J_i$ is  a function of $X^L_i, \xi_i,X^R_{i+1}$, and then from the recursions for $X^L_i$ and $X^R_i$
\[ Z^L_0 = \psi^3 (\xi_K,\xi_{K+1}, \xi_{K+2},\ldots,\xi_{0};  X^L_{K+1}, X^R_1). \]
By (\ref{XxX}) with $X^L_i = X^L_{K+1}$ this is of the form 
\[Z^L_{0} = \psi(\dots,\xi_{-2},\xi_{-1},\xi_0,X^R_{1}) \]
Now (\ref{psi-def}) defines a stationary version of the quadruple process. 
\qed

Just as $X^R_{n,i}$ was the ``looking right" analog of the ``looking left" process $X^L_{n,i}$, 
we can define a ``looking right" process $Z^R_{n,i}$ analogous to $Z^L_{n,i}$ as follows.
Define
 \begin{eqnarray}
 \bV^R_{n,i} = \max_{i \in A \subseteq \{i,i+1,\ldots,n\}} 
\left( |A| - \sum_{j=i}^{n-1} \xi_j \ind(j \in A, j+1 \in A) 
+ \theta \left| (A \bigtriangleup \Aopt)\cap \{i,i+1,\ldots, n \} \right| \right) 
\label{VRn2}
\\
 \bW^R_{n,i} = \max_{i \not\in A \subseteq \{i,i+1,\ldots,n\}} 
\left( |A| - \sum_{j=i}^{n-1} \xi_j \ind(j \in A, j+1 \in A) 
+ \theta \left| (A \bigtriangleup \Aopt)\cap \{i,i+1,\ldots,n \} \right| \right)  .
\label{WRn2}
\end{eqnarray}
Then the difference 
$Z^R_{n,i} = \bV^R_{n,i} - \bW^R_{n,i}$ satisfies the recursion
$$
Z^R_{n,i} = 1 - \min(Z^R_{n,i+1},\xi_i) \ind( Z^R_{n,i+1} >0) + \theta J_{n,i}; \quad 
Z^R_{n,n} = 1 + \theta J_{n,n}.$$ 
Recall that $\Bopt_n$ attains 
$
 \max_{A \subseteq \{1,\cdots,n\}}
\left( |A| - \sum_{i= 1}^{n-1} \xi_i \ind(i \in A, i+1 \in A) + \theta |A \bigtriangleup \Aopt_n|  \right)
$.
As in section \ref{subsec:dp},  
we can write down the benefits of each of the four possible choices for including/excluding items $i$ and $i+1$, and thereby obtain criteria for which combination is used in $\Bopt_n$.  
See Table 3, in which 
$ (Z^L_{n,i},\xi_i,Z^R_{n,i+1})$
is abbreviated to 
$(Z^L,\xi,Z^R)$ and the $n$ subscript is dropped.

\vspace{0.2in} \hspace{-13pt}
{\small
$\begin{array}{llll}
- i - (i+1) - &
\mbox{absolute benefit}&
\mbox{relative benefit}&\mbox{when used} \\
&&\\
- \bullet - - \bullet - &
\bV^L + \bV^R - \xi 
+\theta (\ind_{i \not\in \Aopt}+\ind_{i+1 \not\in \Aopt})
&
Z^L + Z^R - \xi  

& 
\xi <  \min(Z^L - \theta J_i , Z^R - \theta J_{i+1})
\\
&&\quad - \theta (J_i+J_{i+1}) &\\
- \bullet - - \circ - &
\bV^L + \bW^R 
+\theta (\ind_{i \not\in \Aopt}+ \ind_{i+1 \in \Aopt})
&
Z^L 
- \theta J_i 
& 
(Z^R - \theta J_{i+1} )^+ < \min ( Z^L - \theta J_{i}, \xi)
\\
- \circ - - \bullet - &
\bW^L + \bV^R 
+\theta (\ind_{i \in \Aopt}+ \ind_{i+1 \not\in \Aopt})
&
Z^R 
- \theta J_{i+1} 
& 
(Z^L - \theta J_i )^+ < \min ( Z^R - \theta J_{i+1}, \xi)
\\
- \circ - - \circ - &
\bW^L + \bW^R 
+\theta (\ind_{i \in \Aopt}+\ind_{i+1 \in \Aopt})
&
0
& 
\hbox{otherwise}
\end{array}$}

\vspace{0.1in}
\centerline{{\bf Table 3.}  Inclusion criteria for $i, i+1$ in $\Bnopt$.
}

\vspace{0.1in}
\noindent
It should now be clear that the stationary quadruple process can be extended to a 
{\em stationary quintuple process}
$$(Z^L_i , X^L_{i},\xi_i,X^R_{i+1}, Z^R_{i+1}), -\infty < i < \infty  $$ 
in which $Z^R$ satisfies the recursion 
$$Z^R_{i} = 1 - \min(Z^R_{i+1},\xi_i) \ind( Z^R_{i+1} >0) + \theta J_{i}, \ - \infty < i < \infty $$ 
satisfied by $Z^R_{n,i}$.  
By ``reflection symmetry" between $Z^R$ and $Z^L$, 
the functional relationship (\ref{psi-def}) holds for $Z^R$ in reflected form with the same function $\psi$:
%\marginpar{M changed\\ indices\\
%was $\psi(\dots,\xi_{i+2},\\
%\xi_{i+1},\xi_i,X^L_{i-1}$}
\begin{equation} Z^R_{i} = \psi(\dots,\xi_{i+1},\xi_{i},\xi_{i-1},X^L_{i-1})  . \label{psi-def-2}
\end{equation}
We can now use the stationary quintuple process to define  a random subset $\Bopt \subset \Zbold$
by specifing that, for each pair 
$(i,i+1)$, 
we use the one of the four choices which has the largest relative benefit in Table 3.  
Analogously to Lemma \ref{LXinfty} one can check this definition is consistent.
The local weak convergence property 
(Lemma \ref{Lwc})
extends to the present setting as follows.

\iffalse
The optimal set $\Bopt_n$ is obtained by requiring that for each pair $(i,i+1)$, 
we use the one of the four choices which has the largest relative benefit in Table 3, which
translates exactly in the 'when used' rule.
At this stage, it is natural to extend the stationary triple process
 $((X^L_{i},\xi_i,X^R_{i+1}), -\infty < i < \infty )$ which defines $\Aopt$ 
 to a stationary quintuple process $((Z^L_i , X^L_{i},\xi_i,X^R_{i+1}, Z^R_{i+1}), -\infty < i < \infty )$ 
 which define both $\Aopt$ and $\Bopt$ which corresponds to the 'stationary' version of $\Bopt_n$.
Now we define the recursion $$
Z^R_{i} = 1 - \min(Z^R_{i+1},\xi_i) \ind( Z^R_{i+1} >0) + \theta J_{i},$$ 
and we end up with a quintuple process
$$((Z^L_i , X^L_{i},\xi_i,X^R_{i+1}, Z^R_{i+1}), -\infty < i < \infty ). $$
This process has a unique stationary distribution
thanks to Lemma \ref{L4}.  (yyy check lemma)
We can now define consistently a set $\Bopt$
by requiring that, for each pair 
$(i,i+1)$, 
we use the one of the four choices which has the largest relative benefit in Table 3.
%{\tt yyy As in section \ref{sec-fin} we should check this is consistent as $i$ varies -- but I'm sure this will work out.}
\fi

\begin{Lemma}
\label{Lwc2}
Let $U_n$ be uniform on $\{1,\ldots,n\}$.
As $n \to \infty$
\[ 
((Z^L_{n,U_n+i},X^L_{n,U_n+i},\xi_{U_n+i},X^R_{n,U_n+i+1},Z^R_{n,U_n+i+1}, \ind(U_n +i \in \Anopt), \ind(U_n+i \in \Bnopt)), \ - \infty < i < \infty) 
\] \[
\cd 
((Z^L_i,X^L_i,\xi_i,X^R_{i+1}, Z^R_{i+1}, \ind(i \in \Aopt), \ind(i \in \Bopt)), \ - \infty < i < \infty)
 . \]
\end{Lemma}
{\bf Proof.}  The proof repeats the 
proof of Lemma \ref{Lwc}, % using Lemma \ref{L4} 
using (\ref{psi-def},\ref{psi-def-2})
 to incorporate the $(Z^L,Z^R)$ terms.
 %\marginpar{C: we had forgotten to check the continuity of $\ind( i \in \Bopt)$} 
 In order to incorporate the $\Bopt$ component, we need to check that the function $\ind (0 \in \Bopt)$ is a.s. continuous with respect to the stationary distribution of $(Z^L_0,X_0^L, \xi_0,X^R_1, Z^R_1)$. From Table 3, we get that $\{ 0 \in \Bopt \} = \{ Z_0^L - \theta J_0 > \min ( \xi_0 , \max( Z^R_1- \theta J_1 , 0 ) \}$. Hence, it requires that the probability of an equality between some of two $Z^L_0 - \theta J_0, \xi_0, Z^R_1 - \theta J_1$ is zero. We only check that $\Pb ( Z_0^L - \theta J_0 = \xi_0) = 0$. The recursion satisfied by $Z^L_0$ reads $Z^L_{0} - \theta J_{0} = 1 - \min ( Z^L_{-1} , \xi_{-1} ) \ind( Z^L_{-1} > 0)$. Thus, arguing as in the proof of Lemma \ref{L4}, $Z^L_0 - \theta J_0$ is a function of $(Z^L_{K+1}, \xi_{K+1}, \cdots, \xi_{-1}, J_{K+1}, J_{K+2}, \cdots, J_{-1})$ with $K = \max \{ i < 0 : \xi_{i-1} + \xi_i \leq 1 - \tau \}$. Since $Z^L_{K+1} = 1 - \xi_{K} + \theta J_{K+1}$ and $J_i \in \{-1,1\}$, we deduce by recursion that there exists a pair of integers $(i_0,n)$ with $K \leq i_0 \leq  -1$ and $- K \leq n \leq K$ such that $Z^L_0 \in \{ 1 - \xi_{i_0} + n \theta , \xi_{i_0} + n \theta \}$. The independence of $\xi_i$ and $\xi_0$ for $i < 0$ and assumption (\ref{G-assume}) imply that $\Pb ( Z_0^L - \theta J_0 = \xi_0) = 0$.
 \qed

Now define
\begin{eqnarray}
\delta(\theta) \!&\!=\!& \!
\Pb(\{0 \in \Aopt\} \bigtriangleup \{0 \in \Bopt\})
 \label{del-def-0} \\
 \eps(\theta)\!&\!=\!&\!
\Pb(0 \in \Aopt) - \Eb \xi_0 \ind(0 \in \Aopt, 1 \in \Aopt) -\Pb(0 \in \Bopt) + \Eb \xi_0 \ind(0 \in \Bopt, 1 \in \Bopt) .
\label{eps-def-0}
\end{eqnarray}
So $\delta(\theta)$ is the proportion of items at which $\Aopt$ and $\Bopt$ differ, 
and $\eps(\theta)$ is the difference in mean benefit per item between $\Aopt$ and $\Bopt$.
By Lemma \ref{Lwc2}, 
\begin{eqnarray}
\frac{1}{n} \Eb |\Aopt_n\bigtriangleup \Bopt_n|  &=& \Eb |\ind(U_n\in \Aopt_n)-\ind(U_n\in \Bopt_n)|\nonumber \\
&\to & \Pb(\{0 \in \Aopt\} \bigtriangleup \{0 \in \Bopt\}) = \delta(\theta) \label{delta-con}
\end{eqnarray}
and similarly the mean benefits satisfy
\begin{equation}
n^{-1} (\Eb f_n(\Anopt)) - \Eb f_n(\Bnopt)) \to \eps(\theta). \label{eps-con}
 \end{equation}
\begin{Proposition}
\label{le:nMn}
Let  $\bM_n = f_n(  \Bnopt )$ be the benefit associated to $\Bnopt$, then a.s. and in $L^1$
\begin{eqnarray}
\label{eq:nMn1}
\lim_{n \to \infty} n^{-1} | \Bnopt \bigtriangleup \Anopt |  = \delta (\theta) \\
\label{eq:nMn2}
\lim_{n \to \infty} n^{-1} (  M_n -  \bM_n) = \eps (\theta) .
\end{eqnarray}
Moreover for any choice $B^\prime_n$ satisfying (\ref{eq:nMn1}) in $L^1$, the associated benefit $M^\prime_n = f_n ( B^\prime_n)$ satisfies
$$
\liminf_n n^{-1} \Eb (  M_n - M^\prime_n) \geq \eps(\theta).
$$
\end{Proposition}
{\bf Proof.} 
The convergence assertions (\ref{eq:nMn1},\ref{eq:nMn2}) follow from (\ref{delta-con},\ref{eps-con}) and the same concentration argument used in the proof of Theorem \ref{th1};  we will not repeat the details.  By construction, for any $B^\prime_n$ the associated reward $M^\prime_n$ satisfies
$$
M^\prime_n + \theta | B^\prime_n \bigtriangleup \Anopt | \leq \bM_n + \theta | \Bnopt \bigtriangleup \Anopt |.
$$
Then because both $(\Bnopt)$ and $(B^\prime_n)$ satisfy (\ref{eq:nMn1}), we see that
$$
\Eb M^\prime_n \leq \Eb \bM_n + o (n).
$$
\qed

\paragraph{Discussion}
For $0 < \theta < \tau$ and for $\delta = \delta(\theta)$, 
 Proposition \ref{le:nMn} 
implies that the limit 
$\bar{\eps}(\delta) = \lim_n \Eb \eps_n(\delta)$
exists and that 
\[ \bar{\eps}(\delta(\theta)) = \eps(\theta) . \]
So to prove Theorem \ref{th2} it should be enough to prove 
\begin{equation}
\delta(\theta) \sim \alpha \theta, \quad \eps (\theta) \sim \beta \theta^2 
\mbox { as } \theta \to 0
\label{T1}
\end{equation}
for positive constants $\alpha, \beta$.  
Now the definitions (\ref{del-def-0},\ref{eps-def-0}) enable us to rewrite (using Table 3) $\delta(\theta)$ and $\eps(\theta)$ in terms of the stationary distribution  $(Z^L_{0} , X^L_0,\xi_0,X^R_{1},Z^R_{1})$ 
of the quintuple process, as 
\begin{eqnarray}
\delta(\theta) &=&
\Pb \left(
\{X^L_0>\min(X^R_1,\xi_0)\} \bigtriangleup  \{ Z^L_0 - \theta J_0  >  \min (  (Z^R_1-\theta J_1)^+, \xi_0)  \}
\right)   \label{del-def-z} \\
\eps(\theta)&=&
 \Pb(X^L_0 > \min(X^R_1,\xi_0)) - \Pb (Z^L_0 - \theta J_0  >  \min ( (Z^R_1-\theta J_1)^+, \xi_0) ) \nonumber\\
&&  -\  \Eb \xi_0 \left( 
\ind(\xi_0 < \min(X^L_0,X^R_1)) - \ind(\xi_0 < \min (Z^L_0 - \theta J_{0}, Z^R_1 - \theta J_1)  \right)  .  \nonumber
\end{eqnarray}
So if we had an explicit formula for the stationary distribution 
 $(Z^L_{0} , X^L_0,\xi_0,X^R_{1},Z^R_{1})$,
then we could derive an explicit formula for 
$\delta(\theta)$ and $\eps(\theta)$ and
seek to prove (\ref{T1}) by calculus.  
But we do not have such an explicit
formula -- note the independence property (\ref{3-ind}) of the triple process 
does not hold for the quintuple process -- 
and we have not completely succeeded in that program.
We can prove (see Appendix)
the $\delta(\theta) \sim \alpha \theta$ part of (\ref{T1}), 
though we only use the weaker upper bound, proved by a simpler argument 
in section \ref{subsec:existence}
.
To handle $\eps(\theta)$ we show how to
rewrite $\delta(\theta)$ and $\eps(\theta)$ in a different way (Proposition \ref{P-cycle}) that allows us to derive inequalities,
which will establish the stated form
of Theorem \ref{th2}.

\subsection{Existence of the limit function $\bar{\eps}(\delta)$}
\label{subsec:existence}

There is a minor technical point we deal with first.
We expect intuitively that the function $\delta(\theta)$ should be continuous monotone, but neither property is obvious.  
If there were small values of $\delta$ which were not of the form 
$\delta = \delta(\theta)$ for some $\theta$, 
then we can't use Proposition \ref{le:nMn} to establish existence of a limit 
$\bar{\eps}(\delta)$.
Instead we outline an argument 
(reusing previous methods) to prove 
more abstractly
(Lemma \ref{L-lim-exists})
that the limit $\bar{\eps}(\delta)$ always exists.   
We could have started the proof of Theorem \ref{th2} this way, but we wanted to emphasize the Lagrange multiplier approach as more useful for calculation.
\begin{Lemma}
\label{L-lim-exists}
$\bar{\eps}(\delta) : = \lim_{n} \Eb \eps_n(\delta) $ exists, 
for each $0<\delta<1$.
\end{Lemma}
Note that $\eps_n(\delta)$ is {\em a priori} non-decreasing in $\delta$, and hence 
$\bar{\eps}(\cdot)$ is non-decreasing.

\noindent
{\bf Outline proof.}
Fix  $0<\delta<1$.
Let $B_n^{(\delta)}$ attain the minimum in the definition (\ref{eps-def}) of $\eps_n(\delta)$. 
Set $\bar{\eps}_*(\delta) = \liminf_n \Eb \eps_n(\delta)$.  
There exists a subsequence 
(of the subsequence of $n$ attaining the {\em liminf}) 
in which the local weak convergence (Lemma \ref{Lwc}) of $\Anopt$ to $\Aopt$ extends to joint convergence of $B_n^{(\delta)}$ to some limit random set $B^{(\delta)}$.  
The analogs of (\ref{del-def-0}, \ref{eps-def-0}) with $B^{(\delta)}$ in place of $\Bopt$ equal $\delta$ and $\bar{\eps}_*(\delta)$.  
For arbitrary $n$, start with the restriction ($B_n^*$, say) of $B^{(\delta)}$ to $[1,n]$ and then show that by modifying $B_n^*$ near the endpoints we can construct $B_n^{**}$ satisfying 
$|B_n^{**} \bigtriangleup \Anopt| \geq \delta n$ and 
$ \Eb [n^{-1} (f_n(\Aopt_n) - f_n(B_n^{**}))] \to \bar{\eps}_*(\delta)$.\qed

%\marginpar{M: end of section\\
 %changed}
The following Lemma (to be proved in Section \ref{sec:int}) allows to
complete the proof of Theorem \ref{th2}:
\begin{Lemma}\label{lem:star}
There exist positive constants $C_1, C_2$ such that, for all $0 < \theta < \tau$,
\begin{eqnarray}
 \delta(\theta) & \leq &  C_1 \theta 
\label{star1}\\
 \eps(\theta) &  \geq & C_2 \theta^2. 
\label{star2}
\end{eqnarray}
\end{Lemma}
We now finish the proof of Theorem \ref{th2}.
Recall that 
Proposition \ref{le:nMn} showed
$\bar{\eps}(\delta(\theta)) = \eps(\theta)$, 
and that 
(Lemma \ref{L-lim-exists})
$\bar{\eps}(\cdot)$ is a non-decreasing function.
Using (\ref{star1})
\[ \bar{\eps}(C_1 \theta) \geq \bar{\eps}(\delta(\theta))=\eps(\theta)
\geq C_2 \theta^2 \]
and setting $\delta = C_1 \theta$ gives 
$\bar{\eps}(\delta) \geq C_2 \delta^2/C_1^2$.
This establishes 
 the lower bound (\ref{lower-epsilon}) and completes the proof of Theorem \ref{th2}.

\subsection{A cycle formula representation}

\begin{Lemma}
\label{Litau}
If $\xi_{i-1} + \xi_i < 1 - \tau$ then $i \in \Aopt$ and $i \in \Bopt$.
\end{Lemma}
{\bf Proof.}  
Suppose $\xi_{i-1} + \xi_i < 1 - \tau$.  
Lemma \ref{lem1}(a) showed $i \in \Aopt$.  
Recall that $\Bnopt$ maximizes (\ref{iop2}).  
If $i \not\in A$ then the increase in the benefit at (\ref{iop2}) obtained by including $i$ is at least
$1 - \xi_i - \xi_{i-1} - \theta$, so by our standing assumption (\ref{theta-tau}) the increase is positive and so 
$i \in \Bnopt$. 
Letting $n \to \infty$ and using Lemma \ref{Lwc2} gives the same conclusion for $\Bopt$.
\qed

We next need a lemma (analogous to Lemma \ref{lem1}(b)) giving conditions under which we can ``localize" $\Aopt$ and $\Bopt$ by forcing them to coincide with the optimal sets 
$\Anopt$ and $\Bnopt$ for the optimization problem on $[1,n]$ for suitable $n$, which we now write as $t-1$.
\begin{Lemma}
\label{Lrestrict}
Let $t \geq 2$.  Suppose
$\xi_{i-1} + \xi_i < 1 - \tau$ for each of $i = 0,1, t-1, t$.
Then: \\
(a)  $\Aopt$ and $\Bopt$ contain $\{0, 1,  t-1, t\}$. \\
(b) The restrictions of $\Aopt$ and $\Bopt$ to $[1,t-1]$ coincide with $\Atopt$ and $\Btopt$. \\
(c) For  any 
$B \subseteq \{1,2,\ldots,t-1\}$, 
either $B = \Atopt$ or $f_{t-1}(B) < f_{t-1}(\Atopt)$.\\
(d) In particular, either $\Atopt = \Btopt$ or $f_{t-1}(\Atopt) > f_{t-1}(\Btopt)$.
\end{Lemma}
{\bf Proof.}  
(a) follows from Lemma \ref{Litau}.  
Observe that $\Atopt$ and $\Btopt$ contain $1$ and $t-1$, because $\xi_1 < 1 - \tau$ 
and $\xi_{t-2} < 1 - \tau$.  
If we consider the solutions 
$A_{[\ell,m]}^{\mbox{{\tiny opt}}}, \ 
 B_{[\ell,m]}^{\mbox{{\tiny opt}}}$ 
 for some interval $[\ell,m]$ strictly containing $[0,t]$, 
 then they contain $1$ and $t-1$ by the argument for Lemma \ref{Litau}.   
 Thus by optimality the restrictions of $A_{[\ell,m]}^{\mbox{{\tiny opt}}}$ and 
 $B_{[\ell,m]}^{\mbox{{\tiny opt}}}$ to $[1,t-1]$ must coincide with $\Atopt$ and $\Btopt$.  
 So (b) follows from weak convergence, Lemma \ref{Lwc2}.
And (c) follows from the uniqueness result, Lemma \ref{L-unique}.
 \qed

We start by quoting a standard form (cf. \cite{dur91v3} Exercise 6.3.4) of Kac's identity for stationary processes.
\begin{Lemma}
\label{LKac}
Let $(\Xi_i, - \infty < i < \infty)$ be a stationary ergodic sequence on some state space, 
let $\Pb (\Xi_1 \in \bar{D}) > 0$ and let $h(\Xi_1)$ be real-valued and integrable.  
For any $t_0 \geq 1$, define 
$ T = t_0 \min \{i \geq 2: \Xi_{it_0} \in \bar{D}) $.  
Then 
\[ \Eb h(\Xi_1) = \Eb \left[ \ind (\Xi_1 \in \bar{D}) \sum_{i = 1}^{T-1} h(\Xi_i) \right]  .\]
\end{Lemma}
We apply this to 
$\Xi_i = (Z^L_i,X^L_i,\xi_i,\xi_{i-1},\xi_{i-2},X^R_{i+1}, Z^R_{i+1})$
and $t_0 = 3$ and
\begin{equation}
 D := \{\xi_{-1} + \xi_0 < 1 - \tau, \ \xi_0 + \xi_1 < 1 - \tau\} = \{\Xi_1 \in \bar{D}\} 
 \label{D-def}
 \end{equation}
for suitable $\bar{D}$, making the $T$ in Lemma \ref{LKac} be
\begin{equation}
 T =  3\min \{t \geq 2: \xi_{3t-2} + \xi_{3t-1} < 1 - \tau, \ \xi_{3t-1} + \xi_{3t} < 1 - \tau\} . \label{T-def}
 \end{equation}
Now definition (\ref{del-def-0}) says 
$\delta(\theta) = \Eb h(\Xi_0)$ for 
\[ h(\Xi_0) = \ind (\{0 \in \Aopt\} \bigtriangleup \{0 \in \Bopt\}) . \] 
So $\sum_{i = 1}^{T-1} h(\Xi_i) $ equals the cardinality  of 
$ \Aopt \bigtriangleup \Bopt$ restricted to $[1,T-1]$.  
On the event $D$, Lemma \ref{Lrestrict} identifies this restriction as $ \ATopt \bigtriangleup \BTopt$, 
so Kac's identity gives (\ref{delta-cycle}) below.  
Similarly,  definition (\ref{eps-def-0}) says 
$\eps(\theta) = \Eb h(\Xi_0)$ for 
\[ h(\Xi_0) = \ind(0 \in \Aopt) -  \xi_0 \ind(0 \in \Aopt, 1 \in \Aopt) - \ind (0 \in \Bopt) +  \xi_0 \ind(0 \in \Bopt, 1 \in \Bopt)  \] 
and on the event $D$ the sum $\sum_{i = 1}^{T-1} h(\Xi_i) $ equals the difference 
$f_{T-1}(\ATopt ) - f_{T-1}(\BTopt )$ between the benefits.
This establishes (\ref{eps-cycle}), 
and the final assertion (\ref{Dff}) follows from  Lemma \ref{Lrestrict}(d).
To summarize:
\begin{Proposition}
\label{P-cycle}
Let $D$ be the event (\ref{D-def}) and let $T$ be the random time (\ref{T-def}).  Then
\begin{eqnarray}
\delta(\theta) &=& \Eb [ \ind_D \times |\ATopt \bigtriangleup \BTopt| ] \label{delta-cycle}\\
\eps(\theta) &=& \Eb [ \ind_D \times  (f_{T-1}(\ATopt ) - f_{T-1}(\BTopt )) ] . \label{eps-cycle}
\end{eqnarray}
\begin{equation} \mbox{
On $D$, either $\ATopt = \BTopt$ or $f_{T-1}(\ATopt ) - f_{T-1}(\BTopt ) > 0$.
} \label{Dff}
\end{equation}
\end{Proposition}

\subsection{An integration lemma}\label{sec:int}
Let us rewrite the difference in (\ref{eps-cycle}) as 
\[ W(\theta) : = f_{T-1}(\ATopt ) - f_{T-1}(\BTopt ) \]
to emphasize its dependence on $\theta$; and note $D$ does not depend on $\theta$.  
The key ingredient in the proof of the lower bound is the following lemma, to be proved in  section \ref{sec:mu}.

\begin{Lemma}
\label{le:mu}
There exists $  C_3 >0$ such that for all $0 < \theta < \tau$, for all
$k\geq 0$ and $x >0$, 
$$ \Pb ( T\geq k,\:0 < \ind_D W (\theta) < x )  \leq C_3 x(k+1)\Pb(T\geq k) . $$
\end{Lemma}
Taking $k=0$ in this lemma we get
\begin{eqnarray}
\label{eq:linx}\Pb ( 0 < \ind_D W (\theta) < x )  \leq C_3 x.
\end{eqnarray}

Recall a simple integration lemma (for a more general result see \cite{me115} Lemma 6(a)):
\begin{Lemma}
\label{LxiW}
Let $V \geq 0$ be a real-valued random variable % which satisfies for all $x >0$,  $ \Pb ( 0<V <x) \leq C  x$, Then for any event $A \subseteq \{V >
such that
\[ \Pb ( 0<V <x) \leq C  x, \quad 0 < x < \infty . \]
Then
\[ \Eb V \geq \frac{[\Pb (V>0)]^2}{2C} . \]
\end{Lemma}
%Lemma \ref{le:mu} showed 
%$ \Pb( 0 < W (\theta) < x |D) \leq C_1x/\Pb(D)$, 
By (\ref{eq:linx}) and Lemma \ref{LxiW}, we get
\begin{eqnarray}
\eps(\theta) &=& \Eb (\ind_D W(\theta)) \mbox{ by (\ref{eps-cycle})}\nonumber \\
%&=& \Pb(D) \Eb (W(\theta)|D)\nonumber \\
&\geq& \frac{[\Pb( W(\theta)\ind_D > 0)]^2}{2C_3}. \label{eq:epsbound}
\end{eqnarray}
To finish the proof of (\ref{star2}), we need the following lemma 
\begin{Lemma}
\label{lem:contd}
There exists a positive constant $C_4$ such that, for all $0 < \theta < \tau$,
\begin{eqnarray}
\Pb ( W (\theta)\ind_D   > 0 ) \geq C_4 \theta. 
\label{star-2}
\end{eqnarray}
\end{Lemma}
{\bf Proof.}
%\marginpar{M: small modif}
By assumption (\ref{G-assume}) we may assume that the constant $\tau$ 
at (\ref{eq:deftau}) is such that 
\begin{equation}
\inf_{ 1/2 - 2 \tau < x < 1/2  \tau} 
g(x) > 0 
\label{gx-0}
\end{equation}
where $g$ is the density function for $\xi_i$.
Consider the event:
\begin{eqnarray*}
\Omega ( \theta) &  = & \left\{ \xi_{-1} \in (0,1/2) , \; \xi_{0} \in ( 0, 1/2 - \tau), \; \xi_{1} \in ( 1/2  - \tau , 1/2) ,  \right.\\ 
& & \quad \quad \quad \left.\; \xi_{2} \in (  1 - \xi_1 - \theta, 1 - \xi_1) , \;\xi_3 \in ( 0 , 1/2 - 2 \tau ) \right\}
 . \end{eqnarray*}
Using (\ref{gx-0}) there exists $C_4 > 0$ such that
$$
\Pb ( \Omega (\theta)) \geq C_4 \theta.
$$
Assume this event $\Omega(\theta)$ happens.  
Then $\xi_{-1} + \xi_{0} \leq 1 - \tau$, $\xi_{0} + \xi_1 \leq 1 -
\tau$, $1-\theta<\xi_{1} + \xi_{2} < 1$ and $\xi_{2} + \xi_{3} \leq 1
- \tau$.
So $D$ happens and using Lemma \ref{lem1}(a), we have $\{1,2,3\}\in
\Aopt$ and by Lemma \ref{Lrestrict}(b) the same holds true for $\ATopt$.
Still assuming $\Omega(\theta)$ occurs, we see that for $B =\ATopt
\backslash\{2\}$, we have
$f_{T-1}( \ATopt ) - f_{T-1} (B) = 1 - \xi_1 - \xi_2 \in (0 , \theta)$
%and therefore $\Aopt_3 = A$ and $\Bopt_3 = B$,
and therefore $f_{T-1} (B)+\theta|\ATopt \bigtriangleup B|>f_{T-1}( \ATopt )$,
implying 
$0 < W(\theta )$ by (\ref{Dff}).  % It follows that:
In particular
\[ \Pb (W(\theta)\ind_D > 0) \geq \Pb (\Omega(\theta)) \geq C_4 \theta \]
and we have proved the assertion (\ref{star-2}).
\qed

From (\ref{eq:epsbound}) and (\ref{star-2}), we directly get the
second assertion (\ref{star2}) of Lemma \ref{lem:star}. We now show how to obtain the first assertion of Lemma \ref{lem:star}.
Recall that by definition, we have
\begin{eqnarray*}
f_{T-1}(\BTopt) + \theta |\ATopt \bigtriangleup \BTopt|\geq f_{T-1}(\ATopt ),
\end{eqnarray*}
hence we get $\theta T> \theta |\ATopt \bigtriangleup \BTopt|\geq
W(\theta)$. In particular, by Proposition \ref{P-cycle}, we have
$D\cap\{ W(\theta)>0\}\subset D\cap\{ \theta T> W(\theta)>0\}$.
Also by Lemma \ref{le:mu}, we have
\begin{eqnarray*}
\delta(\theta) &\leq& \Eb[T\ind_D\ind(W(\theta)>0)]\mbox{ by (\ref{delta-cycle})}\\
&\leq& \sum_j j\Pb(T\geq j, \theta j > W(\theta)>0)\\
&\leq& C_3 \theta\sum_j j^2(j+1)\Pb(T\geq j)\\
&\leq& C_3 \theta \Eb[(T+1)^4],
\end{eqnarray*}
and $T/3$ has a geometric distribution so that assertion (\ref{star1})
of Lemma \ref{lem:star} follows.

\subsection{Proof of Lemma \ref{le:mu}}
\label{sec:mu}
Write $W = W(\theta)$.
Consider the random collection 
\[ \Bcal(T-1) :=\{B \subseteq \{1,2,\ldots,T-1\}: \  B\neq  \ATopt ,1 \in B, T-1 \in B\} . \] 
By Proposition \ref{P-cycle}
\begin{eqnarray}
\label{ineq:W}
\mbox{on $D$, either  $W = 0$ or }
W\geq \min_{B\in \Bcal(T-1)} 
( f_{T-1}(\ATopt) - f_{T-1}(B)) > 0 .
\end{eqnarray}
Our first goal is to derive a lower bound (Lemma \ref{le:Ww}) for the right-hand side of (\ref{ineq:W}) in terms of the $\xi_i$'s.  
Until the end of the proof of Lemma \ref{le:Ww}, we are working on the event $D$.

Let $C=\arg\min_{B\in \Bcal(T-1)} ( f_{T-1}(\ATopt) - f_{T-1}(B))$ 
be the optimal perturbation of $\Aopt$ on $[1,T-1]$.  
For any subinterval
$\Ical = [\ell,m] \subseteq [1,T-1]$ write 
$\Ical_e = [\max(\ell - 1,1), \min(m+1,T-1)]$.   
Decompose  $\Aopt\bigtriangleup C$ as $ \cup_i \Ical_i$ where the $\Ical_i$'s are disjoint maximal intervals of $\Aopt\bigtriangleup C$.
Then
\[ f_{T-1}(\Aopt)- f_{T-1}(C)= \sum_i \left( f_{(\Ical_i)_e}(\Aopt\cap{(\Ical_i)_e})- f_{(\Ical_i)_e}(C\cap{(\Ical_i)_e}) \right)
= \sum_i \left(  f_{T-1}(\ATopt) - f_{T-1}(C_i) \right) \] 
where $C_i = (\ATopt\cap{\Ical_i}^c) \cup (C\cap{(\Ical_i)_e})$.
%\marginpar{M: change def of $C_i$}
This implies that  $\Aopt\bigtriangleup C$ is a {\em single} subinterval $\Ical$ of $[1,T-1]$.

\iffalse
Note that for any $B\in \Bcal$, we have $f\left(\Aopt|_\Ical\right)-f\left(B|_\Ical\right)>0$ because $\Aopt|_\Ical$ is also  almost surely the unique optimizing set of the (finite) optimization problem of $f$ on $\Ical$.
A similar argument gives the following lemma
\begin{Lemma}\label{le:pos}
If $\Aopt\bigtriangleup B$ is a finite non-empty set, then we have almost surely $f(\Aopt) - f(B)>0$.
\end{Lemma}
Let $C$ denote the optimal perturbation of $\Aopt$ on $\Ical$: $C=\arg\min_{B\in \Bcal}
f\left(\Aopt|_\Ical\right)-f\left(B|_\Ical\right)$. We write $\Aopt\bigtriangleup C= \cup_i \bigtriangleup_i$ where the $\bigtriangleup_i$'s are disjoint maximal intervals, then we have with $C_i=C|_{\bigtriangleup_i}$,
\begin{eqnarray*}
f(\Aopt)- f(C)= \sum_i f(\Aopt)- f(C_i)\leq \min_i f(\Aopt)- f(C_i),
\end{eqnarray*}
and we see that $\Aopt\bigtriangleup C=\bigtriangleup$ is a subinterval of $\Ical$.
\fi

We now look at the possible perturbations of $\Aopt$ on the interval % $\Ical$. We need to introduce some notations first.
$[0,T]$.
Recall that we are working on the event $D$, 
and that $\Aopt$ contains $0,1,T-1,T$.
Let 
$L_0, L_1, \ldots, L_K$ 
be the maximal subintervals 
$[a,b] \subseteq \Aopt \cap [0,T]$ 
for which $b>a$,
that is with at least two elements.  
So we can partition
$[0,T]$ as
$L_0 \cup 
S_0\cup L_1\cup S_1\cup \dots\cup L_K$
where the $S_k$ are the complementary intervals.
We call the $L_k$ 
{\em lakes} 
and we call the $S_k$ 
{\em switches}.

\iffalse
Let $\alpha_0=\inf\{i\geq 1,\: i\notin \Aopt\}$ and $\beta_0=\inf\{i\geq
\alpha_0,\: i+1,i+2\in \Aopt\}$. We call the set
$S_0:=[\alpha_0,\beta_0]$ a switch:
\begin{eqnarray*}
\begin{array}{cccccccccccc}
\alpha_0-1&-&\alpha_0&-&\alpha_0+1&-&&\beta_0&-&\beta_0+1&-&\beta_0+2\\
\bullet&-&\circ&-&\bullet&-&\dots&\circ&-&\bullet&-&\bullet\\
\end{array}
\end{eqnarray*}
Let $\alpha_k=\inf\{i\geq \beta_{k-1},\: i\notin \Aopt\}$ and
$\beta_k=\inf\{i\geq \alpha_k,\: i+1,i+2\in \Aopt\}$. We call the
set $L_k=[\beta_{k-1}+1,\alpha_k-1]$ a lake. We have $\Ical\cup\{I+1\} = L_0\cup
S_0\cup L_1\cup S_1\cup \dots\cup L_K$.
\fi

\begin{Lemma}\label{lem:lake}
Let $L=[a,b]$ be a lake. For any set $B \in \Bcal(T-1)$ such that
$B\cap{L^c}=\Aopt\cap{L^c}$ and 
hence $B\cap{L}\neq \Aopt\cap{L}$, we have
\begin{eqnarray}
\label{ineq:lake}f_{T-1}(\Aopt)-f_{T-1}(B)\geq \min\left\{1-\xi_a, 1-\xi_{b-1},\min_{a\leq i\leq b-1}1-\xi_{i-1}-\xi_i\right\}>0.
\end{eqnarray}
\end{Lemma}
{\bf Proof.}
First suppose $B$ is obtained by removing from $\Aopt$ a single item.  
If this item is $a$, we have $f_{T-1}(\Aopt)-f_{T-1}(B) = 1-\xi_a$; if it is $b$,  we have $f_{T-1}(\Aopt)-f_{T-1}(B)
= 1-\xi_{b-1}$ and if it is $i\in(a,b)$ then we have $f_{T-1}(\Aopt)-f_{T-1}(B) = 1-\xi_{i-1}-\xi_i$. 
So by optimality of $\ATopt$ the first inequality in (\ref{ineq:lake}) holds for these $B$, and Lemma \ref{Lrestrict} implies 
the last inequality in (\ref{ineq:lake}).  
Now recall that 
Lemma \ref{lem1}(c) shows $\min_{a\leq i\leq b-1}1-\xi_i \geq 0$.
So construct a general $B$ by removing items from $\Aopt$ one by one, and for items after the first the benefit can only decrease.
So the first inequality holds generally.  
\qed

\iffalse
Thanks to Lemma \ref{le:pos}, we have $\inf_B f(\Aopt)-f(B)>0$, where the infimum is taken over set $B$ satisfying the conditions of the lemma. Moreover, if $B$ is obtained by removing from $\Aopt$ a single item, then if this item is $a$, we have $f(\Aopt)-f(B) = 1-\xi_a$, if it is $b$,  we have $f(\Aopt)-f(B) = 1-\xi_b$ and if it is $i\in(a,b)$ then we have $f(\Aopt)-f(B) = 1-\xi_{i-1}-\xi_i$. Hence the last inequality in (\ref{ineq:lake}) follows. 
For a general set $B$, first recall that by the first point of Lemma \ref{lem1} $\min_{a\leq i\leq b-1}1-\xi_i>0$ and then construct iteratively the set $B$ by removing points one after another, it is then easy to see that the lemma holds. 
\fi

\begin{Lemma}\label{lem:switch}
Let $S=[a,b]$ be a switch and $S_e=[a-1,b+1]$. For any set $B$ such that
$B\cap{S_e^c}=\Aopt\cap{S_e^c}$ and $B\cap{S}\neq \Aopt\cap{S}$,  we
have
\begin{eqnarray*}
f_{S_e}(\Aopt)-f_{S_e}(B)\geq \min\left\{\min_{a-1\leq i< j\leq b}\xi_i+\xi_j-1,
\min_{a-1\leq i\leq b}\xi_i,\min_{a\leq i\leq b} \xi_{i}-\xi_{a-2},
 \min_{a\leq i\leq b} \xi_{i}-\xi_{b+1}, 1-\xi_{a-2} -\xi_{b+1}\right\}.
\end{eqnarray*}
\end{Lemma}
{\bf Proof.} 
By construction a switch starts and ends with items not in $\Aopt$, and the two items before and after the switch are in $\Aopt$.  
Moreover, Table 2 shows that two adjacent items cannot both be not in $\Aopt$, so the items in a switch 
$[a,b]$ must strictly alternate between in and not-in $\Aopt$, as illustrated in Figure 4.

We first consider a set $B$ obtained from $\Aopt$ by flipping all items in some subinterval 
$[u,v]$ of $[a,b]$.  
There are four  cases, corresponding to whether the endpoints $u,v$ are in or not-in $\Aopt$.
We exhibit these cases below, labeled as e.g. $[\bullet \ldots \circ]$, together with the value of the benefit change 
$f_{S_e}(\Aopt) - f_{S_e}(B)$.

 \begin{eqnarray*}
\begin{array}{cccccccccccccccccccccccc}
&&&&&a&&&&u&&&&&&v&&&&b&&\\
\Aopt&\bullet&-&\bullet&-&\circ&\dots&\bullet&-&\circ&-&\bullet&\dots&\bullet&-&\circ&-&\bullet&\dots&\circ&-&\bullet&-&\bullet\\
B&\bullet&-&\bullet&-&\circ&\dots&\bullet&-&\bullet&-&\circ&\dots&\circ&-&\bullet&-&\bullet&\dots&\circ&-&\bullet&-&\bullet
\end{array}
\end{eqnarray*}

{\bf Figure 4.}
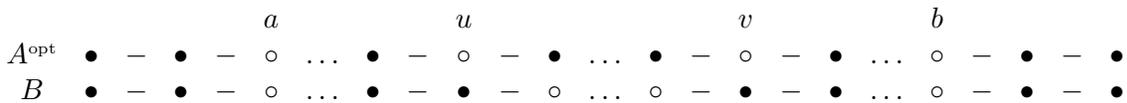
Case $[\circ \ldots \circ]$, where 
$a\leq u\leq v\leq b$.  
Benefit change $=
 \xi_{u-1}+\xi_v-1$.

\begin{eqnarray*}
\begin{array}{cccccccccccccccccccccccc}
&&&&&a&&&&u&&&&&&&&v&&&&b\\
\Aopt&\bullet&-&\bullet&-&\circ&\dots&\bullet&-&\circ&-&\bullet&\dots&\bullet&-&\circ&-&\bullet&-&\circ&\dots&\circ\\
C&\bullet&-&\bullet&-&\circ&\dots&\bullet&-&\bullet&-&\circ&\dots&\circ&-&\bullet&-&\circ&-&\circ&\dots&\circ
\end{array}
\end{eqnarray*}

{\bf Figure 5.}
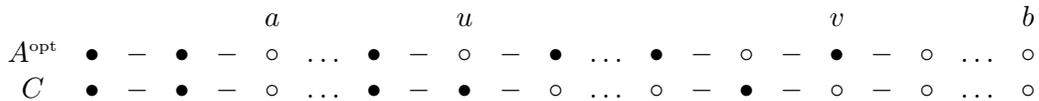
Case $[\circ \ldots \bullet]$, where 
$a\leq u\leq v \leq b-1$.
Benefit change $=
  \xi_{u-1}$.

\newpage
\begin{eqnarray*}
\begin{array}{ccccccccccccccccccccccccccc}
&&&&&a&&&&u&&&&&&&&v&&&&b\\
\Aopt&\bullet&-&\bullet&-&\circ&\dots&\circ&-&\bullet&-&\circ&-&\bullet&\dots&\bullet&-&\circ&-&\bullet&\dots&\circ\\
D&\bullet&-&\bullet&-&\circ&\dots&\circ&-&\circ&-&\bullet&-&\circ&\dots&\circ&-&\bullet&-&\bullet&\dots&\circ
\end{array}
\end{eqnarray*}

{\bf Figure 6.}
Case $[\bullet \ldots \circ]$, where 
$a+1\leq u\leq v\leq b$.
Benefit change $=
 \xi_{v}$.

\vspace{0.08in}
\noindent
In the fourth case $[\bullet \ldots \bullet]$, the benefit change equals $1$.

We also need to consider cases where the flipped subinterval 
$[u,v]$ has $u = a-1$ or $v = b+1$ or both.  
There are five cases, indicated in Table 4.

\vspace{0.08in}

$\begin{array}{llllll}
[u,v] & [a-1, \circ] & [a-1, \bullet] & [\circ, b+1] & [\bullet, b+1] & [a-1, b+1]\\
\mbox{benefit change} & \xi_v - \xi_{a-2} & 1 - \xi_{a-2} & \xi_u - \xi_{b+1} & 1 - \xi_{b+1} & 1 - \xi_{a-2} - \xi_{b+1}
\end{array}$

\vspace{0.05in}
\centerline{{\bf Table 4.} 
}

Now consider any subset $B$ satisfying the hypothesis of Lemma \ref{lem:switch}.   
Decompose $\Aopt \bigtriangleup B$ into disjoint maximal intervals $\Ical_i$.  
It is easy to check that the benefit change between $\Aopt$ and $B$ is just the sum of the separate benefit changes between $\Aopt$ and $\Aopt$ with interval $\Ical_i$ flipped.  
Thus the minimum over $B$ is attained by one of the cases we have considered, 
establishing the lemma.
\qed

\begin{Lemma}
\label{le:Ww}
Set 
$w = \min_{1\leq i< j\leq T-1}\{|\xi_{i}+\xi_j-1|; \xi_i; |1-\xi_i|;
|\xi_i-\xi_j|\}$.  
On the event $D$, either $W = 0$ or $W \geq w$.
\end{Lemma}
{\bf Proof.}
We need only consider the case $W>0$.  
Recall that 
$C=\arg\min_{B\in \Bcal(T-1)} ( f_{T-1}(\ATopt) - f_{T-1}(B))$
is such that  $\Aopt\bigtriangleup C$ is a {\em single} subinterval $\Ical$ of $[1,T-1]$. 
%\marginpar{C: changed $\Delta$ to $\Ical$ to avoid confusion with $\bigtriangleup$}  
It is enough to show that  $C$ satisfies the assumptions of Lemmas \ref{lem:lake} (for some lake) or the assumptions of Lemma
\ref{lem:switch} (for some switch), for then the lower bound $w$ follows from the lower bounds in those lemmas. 

We argue by contradiction: if false, then $\Ical$ intersects some lake and some adjacent switch, say 
$L_k$ and $S_k$ (the case of $L_k$ and $S_{k-1}$ is similar).
So there exist  $a<
b< c$ such that $b= \sup L_k$ and $\Ical = [a, c]$.  
Now check

if $c\in
\Aopt$ then $f(B)-f(C) = 1$ for  
$B:=C\cup\{b,b+2,b+4, \dots,c\}\backslash \{b+1,b+3,\dots, c-1\}$

if $c\not\in
\Aopt$ then $f(B)-f(C) = \xi_c$ for  
$B:=C\cup\{b,b+2,b+4,\dots,c-1\}\backslash \{b+1,b+3,\dots, c\}$.
\\ 
Either case  contradicts the optimality of $C$. 
\qed

We may now complete the proof of Lemma \ref{le:mu}.  
The key point is that the bound $w$ in Lemma \ref{le:Ww} does not depend on $\theta$.
From Lemma \ref{le:Ww},
\begin{eqnarray*}
\Pb (T\geq 3k,\:0 < W (\theta)\ind_D < x )  & \leq & \Pb (T\geq 3k,\:D; 0
< w < x) \leq \Pb(T\geq 3k,\: w<x) \\ 
& \leq & \Pb \left(T\geq 3k,\: \min_{1 \leq i < j \leq T-1} | \xi_i + \xi_j -1 | < x \right) +  \Pb \left(T\geq 3k,\: \min_{1 \leq i  \leq T-1}   \xi_i  < x \right) \\
&& +  \Pb \left(T\geq 3k,\: \min_{1 \leq i  \leq T-1} | \xi_i  -1 | < x \right) +  \Pb \left(T\geq 3k,\: \min_{1 \leq i < j \leq T-1} | \xi_i - \xi_j| < x \right). 
\end{eqnarray*}
The $4$ terms on the right hand side are treated similarly: we will
just study the final term, and will prove that there exists $C>0$
independent of $k$ such that,
\begin{eqnarray}
\label{eq:boundx}
\Pb \left( \min_{1 \leq i < j \leq T-1} | \xi_i - \xi_j | < x \right|
T\geq 3k\Big) \leq C (k+1) x.
\end{eqnarray}
The effect of conditioning 
on the event $\{T\geq 3k\}$  
is that each non-overlapping triple 
$(\xi_{3m},\xi_{3m+1},\xi_{3m+2})$ is conditioned to satisfy either
$\{\xi_{3m}+\xi_{3m+1} \geq 1 - \tau\} \cup\{\xi_{3m+1}+\xi_{3m+2}
\geq 1 - \tau\}$ or $\{\xi_{3m}+\xi_{3m+1} < 1 - \tau,\:\xi_{3m+1}+\xi_{3m+2}
< 1 - \tau\}$ (for $m=T)$. 
It follows that, for any $i<j$,
\begin{equation}
\Pb((\xi_i,\xi_j) \in \cdot | T> j ) 
\leq a^{-2} \Pb((\xi_i,\xi_j) \in \cdot ) 
\label{xi-cond}
\end{equation}
where $a = \min\left(\Pb (\{\xi_{0}+\xi_{1} \geq 1 - \tau\}\cup\{ \xi_{1}+\xi_{2} \geq
1 - \tau\}), \Pb (\xi_{0}+\xi_{1} < 1 - \tau, \xi_{1}+\xi_{2} <
1 - \tau)\right)$.  
From assumption (\ref{G-assume}) the density of $\xi_j - \xi_i$ is bounded by some constant $b$, 
and so
\begin{eqnarray*}
\Pb \left( \min_{1 \leq i < j \leq T-1} | \xi_i - \xi_j | < x \right|
T\geq 3k\Big) &\leq& \sum_{i<j} \Pb(| \xi_i - \xi_j | < x, T\geq
j|T\geq 3k)\\
&=&\sum_{i<j} \Pb(| \xi_i - \xi_j | < x| T\geq
\max(j+1,3k))\Pb(T\geq j+1|T\geq 3k)\\
&\leq& ba^{-2} x\sum_{j\geq 3k-1} (j-1)\Pb(T\geq j+1|T\geq 3k)\\
&\leq& ba^{-2} x\sum_{j\geq k}3(j+1)\Pb(T\geq 3j|T\geq 3k)\\
&=& ba^{-2} x\sum_{j\geq k}3(j+1)\Pb(T\geq 3(j-k))\\
&\leq & ba^{-2} x (k\Eb[T]+\Eb[T(T+1)]),
\end{eqnarray*}
where we used the fact that $T/3$ has a geometric distribution.
This concludes the proof of Lemma \ref{le:mu}.

\section{Final remarks}
\subsection{Technical assumptions on $G$}
\label{sec-G-assume} 
We stated a single assumption (\ref{G-assume}) on $G$.  
What we actually used was three consequences of this assumption:
\begin{itemize}
\item $\Pb (\xi < 1/2) > 0$, which implies $\Pb (\xi_i + \xi_{i+1} < 1) > 0$,  used in Lemma \ref{Litau} and thereby
throughout section \ref{sec:lowerbound} (because it implies $i \in \Aopt$) to implement
``localization" arguments.
\item $\Pb (\xi \leq 1/2) < 1$, used in section \ref{sec-pr-ub} to show $P(\Omega_g) > 0$. 
Note that if $\Pb (\xi \leq 1/2) = 1$ then the optimization problem is degenerate in that 
the optimal $\Anopt = \{1,2,\ldots,n\}$.
\item $\xi_1 + \xi_2$ has density bounded below in some interval $(1,1 + \eta)$, 
used in section \ref{sec-pr-ub} to obtain (\ref{qr-lim}).
\end{itemize}
The latter two are used only in a convenient way to exhibit one near-optimal set.  
The ``localization"  arguments essentially just require one to find some event of positive probability involving $(\xi_{-k},\ldots,\xi_k)$ which forces items $0$ and $1$ to be in (or not in) $\Aopt$.  
Lemma \ref{Litau} is just a simple way to exhibit such.  
So we expect Theorem \ref{th2} to remain true
under much weaker assumptions on $G$.

\subsection{Parallels with the cavity method}

\label{sec-cavity}
The unsophisticated arguments in this paper in the context of i.i.d.-DP (dynamic programming) may be compared with the more sophisticated arguments from the statistical physics {\em cavity method} \cite{MParisi03}, as reformulated in more probabilistic language in \cite{me94,me101}, whose prototype example we take to be the analysis of TSP in the ``mean-field" model of geometry where there are $n$ points and each of the 
${n \choose 2}$ inter-point links has random length.  
Of course {\em algorithmically} DP and TSP are quite different, but there are striking parallels between the analysis of optimal solutions of iid-DP and mean-field-TSP, as follows.
\begin{itemize}
\item There are $n \to \infty$ limits for the random data; in DP this is just the obvious infinite i.i.d. sequence, while for mean-field-TSP it is a certain random infinite tree.  
\item The ``inclusion criterion" for iid-DP involves $X^L_i, X^R_{i+1}$ and the edge-cost $\xi_i$.
Finite-$n$ TSP has of course no simple inclusion criteria, but in the $n \to \infty$ limit of mean-field-TSP there is an analogous criterion for inclusion of an edge $(i,j)$ in terms of quantities 
$Z^L_i,Z^R_j$ and the edge-length $\xi_{ij}$.
Each $Z$ is interpreted (cf. (\ref{XRVW}) for DP) as the difference between costs of two optimal solutions (subject to different local constraints) on one side of the tree. 
\item The distribution we use for $X$ in iid-DP, the stationary distribution of a Markov chain, 
is the solution of an equation with  abstract structure 
$X \ed h(\xi, X^1)$.  
The distribution we use for $Z$ in mean-field-TSP, by a recursion on the limit tree, 
is the solution of an equation with  abstract structure 
$Z \ed h(\xi; Z^1, Z^2,Z^3, \ldots)$ 
where the $Z^j$ are i.i.d. copies of the unknown distribution $Z$.
\end{itemize}
These parallels provide a glimpse of how the analog of Theorem \ref{th1}, a formula for the asymptotic expected cost in mean-field-TSP, may be derived (the original non-rigorous argument was in \cite{MP86}; a rigorous proof was given only recently via more combinatorial methods \cite{wast-tsp}). 
The analog of Theorem \ref{th2} for mean-field-TSP, using Lagrange multipliers as in this paper, 
and leading to a non-rigorous argument that the scaling exponent equals $3$, was given in \cite{me103}.

\paragraph{Acknowledgement}
We thank Vlada Limic for discussions regarding the NK model.

\newpage
%\bibliographystyle{plain}
%\bibliography{../../trees/alg,../../trees/networks,../../trees/me,../../trees/trees,../../trees/misc,../../trees/rwgbook,../../trees/biology}
%\bibliography{/saruman/accounts/fac/pitman/search/general,/saruman/accounts/fac/pitman/search/bm3,/saruman/accounts/fac/pitman/search/bm4,/saruman/accounts/fac/pitman/search/bessel,/saruman/accounts/fac/pitman/search/sizebias,/saruman/accounts/fac/pitman/search/pitman,/saruman/accounts/fac/pitman/search/comb,/saruman/accounts/fac/pitman/search/species,/saruman/accounts/fac/aldous/trees/me,/saruman/accounts/fac/aldous/trees/coag,/saruman/accounts/fac/aldous/trees/trees,/saruman/accounts/fac/aldous/trees/rwgbook,/saruman/accounts/fac/aldous/trees/alg,/saruman/accounts/fac/aldous/trees/misc,/saruman/accounts/fac/aldous/trees/biology}

\newpage
\appendix
\section{Exact asymptotic for the near-optimal solution}
\label{sec:variational}
In this section, we go one step further in the analysis of the quintuple process and prove the following refinement:
\begin{Proposition}
\label{le:delta}
Under the assumptions of Theorem \ref{th2},  $\theta \mapsto \delta (\theta) $ is differentiable at $0$ and $\delta'(0) = \alpha > 0$.
\end{Proposition}
The proof is based on a variational analysis which may have an interest in its own. In view of Proposition \ref{le:delta}, it is tempting to conjecture that there exists $\beta >0$ such that
\begin{equation}
\label{eq:beta}
\eps (\theta) = \beta \theta^2 + o(\theta^2),
\end{equation}
However, we do not have a rigorous proof of this claim. If (\ref{eq:beta}) was true then, from Proposition \ref{le:nMn}, $\delta \mapsto \bar\eps (\delta)$ would be twice differentiable at $0$ and $\bar\eps(\delta) = \beta \alpha^{-2} \delta^2 + o(\delta^2)$. We will exhibit formulae for the constants $\alpha$ and $\beta$.  In order to prove Proposition \ref{le:delta}, using (\ref{del-def-z}) we will write a first order  expansion of $\delta(\cdot)$ at $0$.  We define 
$$
S_i^L (\theta) = Z_i^L (\theta) - \theta J_i - X_i^L \quad \hbox{ and } \quad S^R_{i+1} (\theta)= Z^R_{i+1} (\theta) - \theta J_{i+1} - X^R_{i+1}.
$$
In the next lemma, we give a first order estimate of $S_0^L(\theta)$ and $S_{1}^R(\theta)$ as $\theta$ goes to $0$. Recall that the positive number $\tau$ was defined by (\ref{eq:deftau}). Finally let $K^L =   \inf \{ i \geq  2 : \xi_{-i}  +  \xi_{-i+1} < 1 - \tau \}$. Note that $K^L$ is  independent of $(\xi_0,X_1^R)$.
\begin{Lemma}
\label{le:Q}
There exist a constant $C_5$ and an integer valued random variable $Q^L_0$ independent of $\theta$ such that, for all $0 \leq \theta < \tau$,
$$
|S^L_0 (\theta)| \leq \theta K^L \quad \hbox{ and } \quad \Eb\left[ K^L \ind( S^L_0 (\theta) \neq \theta Q_0^L  ) \right]\leq C_5 \theta,
$$
and respectively for $S^R_1 (\theta)$ and $K^R =  \inf \{ i \geq  2 : \xi_{i}  +  \xi_{i-1} < 1 - \tau \}$.
\end{Lemma}
{\bf Proof.}
To simplify notation in the proof, we set $K = K^L$ and $S_i ^L = S^L_i (\theta) $. By Lemma \ref{L4}, $X^L_{-K+2} = 1 - \xi_{-K+1}$ and $S^L_{-K+2} = 0$. 
Now, since the mapping $x \mapsto \min(x,\xi) \ind(x >0)$ is $1$-Lipschitz (i.e. contracting), we get:
\begin{eqnarray*}
\left |S^L_{i+1} \right| = \left| \min( Z^L_{i} , \xi_i)\ind(Z^L_{i} >0) - \min (X^L_i , \xi_i)\ind(X^L_{i} >0) \right|
 \leq  | Z^L_{i} - X^L_i | 
 \leq  \theta + |S^L_i|.
\end{eqnarray*}
In particular: $
|S_0^L | \leq  \theta (K-2) $,
and it concludes the first statement of the lemma. Now, using (\ref{eq:recZ}), we obtain the recursion
\begin{eqnarray*}
S_{i+1}^L = \min(X_i^L,\xi_i) - \min ((X_i^L+ S^L_i + \theta J_i)^+,\xi_i).
\end{eqnarray*}
Therefore,\begin{itemize}
\item[-] if $X^L_i\geq \xi_i$ then $i\in \Aopt$ and $J_i=-1$. Moreover if
  $X^L_i\geq\xi_i+\theta-S^L_i$, then $S^L_{i+1}=0$.
\item[-] if $X^L_i<\xi_i$ and $-S^L_i-\theta J_i\leq X^L_i<
  \xi_i-S^L_i-\theta J_i$, then we have
  $S^L_{i+1}=-S^L_i-\theta J_i$.
\end{itemize}
For each integer $k$, we define the event $E_k $ as 
\begin{eqnarray*}
E_k = \cap_{-k\leq i\leq -1}\left\{ |X^L_i-\xi_i|\geq k\theta, \:X^L_i\geq
k \theta \right\}.
\end{eqnarray*}
In particular, on the event $E_k \cap \{S^L_{-k} = 0 \}$, we have for $-k\leq i\leq -1$,
\begin{eqnarray*}
S^L_{i+1} = \left\{\begin{array}{ll}
0 &\mbox{ if $X^L_i > \xi_i$,}\\
-S^L_i-\theta J_i &\mbox{ otherwise.}
\end{array}\right.
\end{eqnarray*}
We define $
U = \inf\{i\geq 1, \: X^L_{-i} \leq \xi_{-i}\}$. 
Since $S^L_{-K+2}=0$,  on the event $E_{K-2}$, 
$$S^L_0 =\ind(U\leq K-1)\theta \sum_{i=-U}^{-1}(-1)^i J_{i} = \theta Q^L_0.$$
It remains to upper bound the probability of the event $E^c_{K-2}$. Note that we have
$X^L_{-K+2}=1-\xi_{-K+1}$ so that by the recrusion, we get for any
$k\geq 2$,
\begin{eqnarray*}
X^L_{-K+k} \in \{ \xi_{-K+1},\dots, \xi_{-K+k-1}\}\cup \{1-\xi_{-K+1},\dots, 1-\xi_{-K+k-1}\}.
\end{eqnarray*}
This implies the following inclusion:
\begin{eqnarray*}
\underbrace{\bigcap_{-K +1 \leq   j  \leq -1}  \hspace{-10pt} \{ \xi_j  \notin [0,
K \theta] \cup [1-K \theta, 1]  \} \hspace{-10pt}\bigcap_{-K+1  \leq
  j < i \leq -1}  \hspace{-10pt}\{ | \xi_{j} - \xi_i| \geq  K \theta
\}  \cap \{ | 1- \xi_{j} - \xi_i| \geq  K \theta  \}}_{F_{K}}\subseteq E_{K-2}.
\end{eqnarray*}
We have
\begin{eqnarray*}
\Pb(F_{K}^c,\: K=k) &\leq& \sum_{j= -k+1}^ {-1}  \Pb ( \xi_j  \in [0, k
\theta] \cup [1-k \theta, 1],\: K = k ) \\
&& \quad  + \sum_{-k  +1 \leq   j < i \leq -1} \Pb( |\xi_j - \xi_i |
\in [0, k  \theta] \cup [1 - k \theta,1+ k \theta],\: K = k )  
\end{eqnarray*}
Given $K = k$, the variables $(\xi_i, -k+2 \leq i \leq -1)$ are conditioned on the event $\{\min ( \xi_{i} + \xi_{i-1} , \xi_i + \xi_{i+1} ) \geq 1 -\tau \}$. Similarly, the variable $\xi_{-k+1}$ is conditioned on $\{ \xi_{-k+1} + \xi_{-k} < 1- \tau \} \cap \{ \xi_{-k+1} + \xi_{-k+2} \geq 1- \tau\}$ while the variable  $\xi_{-k}$ is conditioned on $\{ \xi_{-k+1} + \xi_{-k} < 1- \tau \}$. Since we condition  on events of positive probability and using the assumption that $\xi$ has a bounded density, we obtain, for some constant $C$, for all $-k +1 \leq j < i  \leq -1$,
\begin{eqnarray*}
\Pb ( \xi_j  \in [0, k
\theta] \cup [1-k \theta, 1],\: K = k ) & \leq & C k  \theta \Pb (K = k),\\
\Pb( |\xi_j - \xi_i |
\in [0, k  \theta] \cup [1 - k \theta,1+ k \theta],\: K = k )  & \leq & C k  \theta \Pb (K = k).
\end{eqnarray*}
Finally, we get,
\begin{eqnarray*}
\Pb( S^L \neq \theta Q^L , K = k )  \leq    \Pb ( F^c_T , \;  K = k )   \leq   C'  k^3 \Pb ( K = k) \theta. 
\end{eqnarray*}
From (\ref{eq:deftau}), $\Eb [K^3] < \infty$, therefore $\Eb[ K \ind( S^L_0 \neq \theta Q^L_0 ]  \leq C' \Eb [ K^3 ]  \theta$.
\qed

Now, as usual let $(X^L,\xi,X^R) := (X^L_0, \xi_0,X^R_1)$ and similarly, we drop the indices of $S^L_0,Q^L_0,S^R_1,Q^R_1$. Let 
$$\Vcal = \{
(x^L,z,x^R) \in \R^3 : x^L > \min ( z,\max(0,x^R)) \},$$ and $\partial \Vcal = \Vcal^c \cap \overline \Vcal$ is the boundary of $\Vcal$. Note that, from Tables 2 and 3,
\begin{eqnarray*}
\label{eq:idee}  \{ 0 \in \Aopt \} & = & \{( X^L ,\xi,X^R) \in \Vcal \},\\
\{ 0 \in \Bopt \} & = & \{( X^L + S^L,\xi,X^R+S^R) \in \Vcal \}.\nonumber \end{eqnarray*}
We define the functions
\begin{eqnarray*}
F(\theta) & = &  \Pb ( ( X^L,\xi,X^R)  \in  \Vcal ,  ( X^L + S^L(\theta),\xi,X^R+S^R(\theta) )  \not\in  \Vcal  ), \\
\overline F(\theta) & = &  \Pb ( ( X^L,\xi,X^R)  \in  \Vcal ,  ( X^L + \theta Q^L,\xi,X^R+\theta Q^R)  \not\in  \Vcal  ).
\end{eqnarray*}
The next lemma states that $F$ and $\overline F$ have the same first order asymptotic as $\theta$ goes to $0$.
\begin{Lemma}
\label{le:Fbar}
As $\theta \downarrow 0$, 
$$
| F(\theta) - \overline F(\theta) | = o (\theta).
$$
\end{Lemma}
{\bf Proof.}
For $\underline x \in \R^3$ and $r>0$, let $B(\underline x, r)$ denote the closed ball of radius $r$ and center $\underline x$. From Lemma \ref{le:Q},
$$| F(\theta) - \overline F(\theta) |  \leq \Pb \left( \theta(Q^L,Q^R) \ne (S^L,S^R) ; B ( (X^L, \xi,X^R) , \theta \max(K^L,K^R) ) \cap \partial \Vcal \neq \emptyset \right).$$   
Using the "reflection symmetry" from $"L"$ to $"R"$ we obtain:
\begin{eqnarray*}
| F(\theta) - \overline F(\theta) | & \leq & 2 \sum_{k=2} ^ {\infty} \Pb \left( \theta Q^L \ne S^L ; K^L = k ; B ( (X^L, \xi,X^R) , \theta \max(k,K^R) ) \cap \partial \Vcal \neq \emptyset \right).
\end{eqnarray*}
Since the boundary of $\Vcal$ is smooth and the variables $(X^L,\xi,X^R)$ are independent with bounded density, (see (\ref{G-assume},\ref{dist-1})), a simple calculation shows that there exists a constant $c>0$ such that for all $x^L \in [0,1]$ and $t >0$,
$$
\Pb \left( B ( (x^L, \xi,X^R) , t ) \cap \partial \Vcal \neq \emptyset \right) \leq c t. 
$$
Now, fix an integer $n$, from the independence of $(\xi,X^R,K^R)$ and $(Q^L,S^L,X^L,K^L)$, we get 
\begin{eqnarray*}
| F(\theta) - \overline F(\theta) | & \leq & 2 \sum_{k=2} ^ {\infty} \Pb \left( \theta Q^L \ne S^L ; K^L = k ; B ( (X^L, \xi,X^R) , \theta \max(k,n) ) \cap \partial \Vcal \neq \emptyset \right) \\
& & \quad \quad \quad + \; \Pb \left( \theta Q^L \ne S^L ; K^R > n \right) \\
& \leq & 2 \sum_{k=2} ^ {\infty} c \theta \max(k,n)  \Pb \left( \theta Q^L \ne S^L ; K^L = k \right) + C_5 \theta \Pb ( K^R >n ) \\
& \leq & 2 C_5 c n \theta^2 + C_5 \theta \Pb ( K^R > n).
\end{eqnarray*}
Taking $n$ arbitrary large, we obtain our result. \qed

We may now conclude the proof of Proposition \ref{le:delta}.

\noindent {\bf Proof of Proposition \ref{le:delta}.}
We have
\begin{eqnarray*}
\nonumber \overline F(\theta) & = & \int f_3(x^L,z,x^R) \ind_{(x^L,z,x^R) \in \Vcal}  \Eb_{(x^L,z,x^R)}  \ind_{(x^L
+\theta  Q^L,z,x^R+\theta Q^R) \notin \Vcal} dx^Ldzdx^R,\\
\label{eq:F}  & = & \int f_3(\ux) \ind_{\ux \in \Vcal}  \Eb_{\ux}  \ind_{\ux + \theta Q \notin \Vcal} d\ux \label{eq:Ftheta}
\end{eqnarray*}
where $\ux = (x^L,z,x^R)$, $Q = (Q^L, 0,Q^R)$, $f_3 (x^L,z,x^R ) = f(x^L) g(z) f(x^R)$ is the density of the triple process, and
$\Eb_{(x^L,z,x^R)}[\, \cdot \,] = \Eb[ \, \cdot \,| (X^L,\xi,X^R) = (x^L,z,x^R)]$.  Let $\overrightarrow{n}$ be the oriented orthonormal vector to the surface $\partial \Vcal$ at $\ux$, and $dS$ the Lebesgue measure on $\partial \Vcal$. Since $\Vcal^c$ is the intersection of hyperplanes and $Q$ is integer valued, calculus gives
$$
\overline F'(0) = \int_{\partial \Vcal} f_3(\ux)   \Eb_{\ux}  [\max( \overrightarrow{n} .Q, 0) ]  dS,
$$
where $\overrightarrow{n} .Q$ is the usual scalar product of $\overrightarrow{n}$ and $Q$. From Lemma \ref{le:Fbar}, we get $F'(0) = \overline F'(0)$. Recall that (\ref{del-def-z}) implies that
$$
\delta(\theta) =  F(\theta) + \Pb ( ( X^L,\xi,X^R)  \notin  \Vcal ,  ( X^L + S^L(\theta),\xi,X^R+S^R(\theta) ) \in  \Vcal  ).
$$
The same computation on the second term gives $\delta (\theta) = \alpha \theta +  o(\theta)$ with
\begin{eqnarray*}
\alpha &  =  & \int_{\partial \Vcal}  f_3(\ux)   \Eb_{\ux}  | \overrightarrow{n} .Q | dS \\
& = & \int_{0}^{1}\int_{x}^{\infty} f_3 (x,z,x) \Eb_{(x,z,x)} |Q^L - Q^R | dz dx + \int_{0}^{1} \int_{x}^{1} f_3 (x,x,x^R) \Eb_{(x,z,x)} |Q^L| dx dx^R 
\end{eqnarray*}
and this concludes the proof of Proposition \ref{le:delta}. \qed

We now discuss the conjecture (\ref{eq:beta}). Let $\eps_1 (\theta) = \Pb( 0 \in \Aopt) - \Pb ( 0 \in \Bopt) = \Pb ( 0 \in \Aopt \backslash \Bopt ) -  \Pb (0 \in \Bopt \backslash \Aopt )$, we get
\begin{eqnarray*}
\eps_\ind(\theta) & = &  \Pb ( ( X^L,\xi,X^R)  \in  \Vcal ,  ( X^L + S^L,\xi,X^R+S^R)  \not\in  \Vcal  ) \\
 && \quad \quad -\Pb ( (
X^L,\xi,X^R) \not\in  \Vcal ,  ( X^L + S^L,\xi,X^R+S^R)  \in  \Vcal  ),
\end{eqnarray*}
and from what precedes, we deduce that $\eps'_\ind(0) = 0$. Now, we define
$$\Wcal = \{
(x^L,z,x^R) \in \R^3 : z <  \min (x^L, x^R)  \},$$
so that, from Tables 2 and 3,\begin{eqnarray*}
\label{eq:idee2}  \{  \{0,1\} \subset \Aopt \} & = & \{( X^L ,\xi,X^R) \in \Wcal\},\\
\{ \{0,1\} \subset \Bopt \} & = & \{( X^L + S^L,\xi,X^R+S^R) \in \Wcal \}.\nonumber \end{eqnarray*}
Similarly, we introduce the function
\begin{eqnarray*}
\eps_2(\theta) & = &  \Eb [  \xi \ind_{  \{0,1\} \subset \Aopt} \ind_{  \{0,1\}
\not\subset \Bopt} ]   -  \Eb [ \xi \ind_{  \{0,1\}\subset \Bopt} \ind_{  \{0,1\} \not\subset \Aopt} ] \\
 & = &  \Eb [  \xi \ind_{( X^L,\xi,X^R)  \in  \Wcal}\ind_{ ( X^L + S^L,\xi,X^R+S^R)  \not\in  \Wcal  } ] \\
 && \quad \quad \quad - \Eb[  \xi \ind_{ (
X^L,\xi,X^R) \not\in  \Wcal}\ind_ { ( X^L + S^L,\xi,X^R+S^R)  \in  \Wcal  }].
\end{eqnarray*}
We might prove again that $\eps_2'(0) = 0$. Note that
$$
\eps (\theta) = \eps_\ind(\theta)  - \eps_2(\theta) .
$$
Let $$
G (\theta) = \Eb  \xi 1 ( X^L,\xi,X^R)  \in  W ,  ( X^L + S^L(\theta),\xi,X^R+S^R(\theta) )  \not\in  \Wcal  ).
$$
If we had proved that $F$ and $G$ have a second derivative at $0$, then we could obtain
$$
\eps(\theta) = (F''(0) - G''(0)) \theta^2 + o(\theta^2) = \beta \theta^2 + o (\theta^2).
$$
However the computation of the second derivatives of $F$ and $G$ involves some technicalities that we will not consider in this paper.

\end{document}